\documentclass[11pt,a4paper]{article}
\usepackage{indentfirst}
\usepackage{fancyhdr}
\usepackage{color}
\usepackage{mathrsfs}
\usepackage{amsfonts,amssymb,amsmath,amsthm}
\usepackage{authblk}
\usepackage[unicode,pdftex]{hyperref}
\hypersetup{
	colorlinks = true,
	citecolor = green,
	anchorcolor = blue,
	linkcolor = blue}

\newtheorem{thm}{Theorem}[section]
\newtheorem{lemma}[thm]{Lemma}
\newtheorem{coro}[thm]{Corollary}
\newtheorem{prop}[thm]{Proposition}

\newtheoremstyle{rem}{10pt}{10pt}{\rmfamily}{}{\bfseries}{.}{.5em}{}
\theoremstyle{rem}
\newtheorem{rem}[thm]{Remark}

\numberwithin{equation}{section} %该行命令需要放在 hyperref 之下, 否则公式跳转 (x.y.z) 是跳向 (x.1.z)

\title{Two types of averaging principles for multiplicative stochastic Korteweg-de Vries equations with fast oscillations}
\author[1]{Jie Chen}
\author[2]{Fan Gu \thanks{Corresponding Author: gufan@amss.ac.cn}}
\author[3]{Boling Guo}
\author[4]{Xianliang Zhao}
\affil[1]{\scriptsize \textit{School of Science, Jimei University, Xiamen 361021, P.R. China}}
\affil[2]{\scriptsize \textit{School of Statistics and Mathematics, Central University of Finance and Economics,  Beijing 102206, P.R. China}}
\affil[3]{\scriptsize \textit{Institute of Applied Physics and Computational Mathematics, Beijing 100088, P.R. China}}
\affil[4]{\scriptsize \textit{
		Beijing International Center for Mathematical Research, Peking
		University, Beijing 100871, P.R. China}}
\date{}

\begin{document}
	\maketitle
\begin{abstract}
In this paper, we study  averaging principles for stochastic KdV equations perturbed by a fast oscillating term generated by a general stochastic dissipative   equation in $H^s(\mathbb{R})$, $ s\geq 0$. When dealing with multilinear estimates of the nonlinear terms in the stochastic KdV type equation,  two kinds of workspace are both effective. The classical workspace helps us prove the averaging principle with a convergence rate in the $L_\omega^2$ sense, but it has more restrictions on noise kernels and  conservation structures.  When we use a new  workspace without homogeneous low frequency restrictions,  only the regularity conditions and the $L_x^2$ conservation structure is needed. Correspondingly, because there is no moment estimate of the solution, we can only prove a weak type averaging principle in the probability sense. 
Our works may have important reference value, because it not only provides a new perspective to handle dynamical  problems of stochastic dispersive models, but also  compare its advantages with the classical methods.

\end{abstract}		
\section{Introduction}\label{sec:introduction}

In this paper, we consider averaging principles of the following two-time-scale stochastic multiplicative equations:
\begin{equation}\label{moskdv}
	\left\{
	\begin{aligned}
		&d u_\delta=-\left(\partial_{xxx}u_{\delta}+\frac{1}{2}\partial_x(u_{\delta}^2)+u_\delta v^{2}_\delta\right)dt +u_\delta\Phi_1 dW_t^1,\\
		&dv_\delta=\delta^{-1}\left(\mathcal{L}v_\delta-v_\delta\right)dt+\delta^{-1/2}v_\delta\Phi_2 dW_t^2;\\
		& (u_{\delta},v_{\delta})|_{t = 0} = (u_0,v_0),
	\end{aligned}
	\right.
\end{equation}
defined on a given complete filtrated probability space $\left(\Omega, \mathcal{F},\mathbb{P},\left\{\mathcal{F}_t\right\}_{t\geq 0} \right)$. This is a continuation of our series of work \cite{kdvas} and \cite{kdv_sharper} on stochastic KdV equations.

Let us explain the notations in \eqref{moskdv}.  $\delta\in(0,1] $ is a small parameter describing the ratio of the time scale between the slow component $u_\delta$ and the fast component $v_\delta$. $\mathcal{L}$ is a general dissipative multiplier, which means $ \mathscr{F}(\mathcal{L}f)(\xi)=L(\xi)\mathscr{F}(f)(\xi)$. Here,  $L(\xi)$ is an arbitrary fixed non-positive measurable function.

For example, $\mathcal{L}v_\delta$ can be $\partial_{xx} v_\delta$. In this case, $v_\delta$ in \eqref{moskdv} is a solution of a stochastic heat type equation. 

According to our well-posedness result in \cite{kdvas}, we let
\begin{equation*}\label{kernel}
	\Phi_i f(x) = \int_{\mathbb{R}} k_i(x,y)f(y)~dy,\ i=1,2
	,\  \forall~f\in L^2(\mathbb{R}).
\end{equation*}

What's more, $W^i_t=\sum_{j=0}^{\infty}\beta^i_j(t)e_j(x),\ i=1,2$, where  $\left\{\beta^i_j\right\}_{j\in\mathbb{N}},i=1,2$ are two sequences of mutually independent standard Brownian motions and  $\left\{e_j\right\}_{j\in\mathbb{N}}$ is an orthonormal basis of $L^2(\mathbb{R})$.

The mild solution of \eqref{moskdv} can be denoted by 
\begin{equation}\label{mild_solution}
	\left\{\begin{aligned}
	u_\delta=& U(t)u_0-\int_{0}^{t} U(t-r) \left(\frac{1}{2}\partial_{x}(u_\delta)^{2} +u_\delta v^{2}_\delta\right)dr+\int_{0}^t U(t-r) (u_\delta\Phi_1 dW_r^1), \\
	v_\delta=&S(t)v_0-\delta^{-1}\int_{0}^{t} S(t-r) v_\delta dr+\delta^{-1/2}\int_{0}^{t} S(t-r) (v_\delta \Phi_2 dW_r^2), 
	\end{aligned}\right.
\end{equation} 
where $U(t)=e^{-t\partial_{xxx}}$ and $S(t)=e^{\frac{t}{\delta}\mathcal{L}}$(Although $S(t)$  depends on $\delta$, we still use the notation $S(t)$ for simplicity).

When we consider multiscale systems arising from mechanics, mathematics, physics,
chemistry, biology and other areas, averaging principles play an important role in investigating multiscale nonlinear dynamical systems. It provides a powerful tool for simplifying multiscale nonlinear dynamical systems and obtaining approximate solutions. Usually, a multiscale model can be described through coupled equations with ``slow"
and ``fast" components. The averaging principles illustrate that the slow component will converge to a so-called averaged equation as the ratio $\delta$ down to zero. 

The averaging principle for deterministic systems was first studied in \cite{bogoliubov61}, which then was extended to stochastic ordinary differential equations in \cite{averagingSDE}. After the averaging principle was established for a general kind of reaction-diffusion with two time-scales in \cite{aprecdiff}, the averaging principle has been studied for a lot of kind of SPDEs. For instance, stochastic reaction-diffusion equations with polynomial coefficients, stochastic wave equations, stochastic Burgers equations, stochastic two dimensional Navier–Stokes equations were studied in \cite{aprecdiffpoly}, \cite{aphpeq}, \cite{apburgers} and \cite{ap2dns}, respectively. There are also works about abstract stochastic equations and dispersive equations defined on bounded domain, such as \cite{apabstract}, \cite{apschrodingerbounded} and \cite{apkdvbounded}.

However, because the  well-posedness results of dispersive equations in $\mathbb{R}$ are very different from the models given in above literatures and references therein (owing to the  non-dissipative property and nonlinear terms), as far as we know, there have been no result about the averaging principle for stochastic dispersive equations defined on $\mathbb{R}$ until now. Thus, in order to study the interaction processes of KdV equations and oscillations systems in random and multiscale environment, we study the averaging principle of \eqref{moskdv}. Furthermore, our work also have important reference value, when others consider the dynamical properties of stochastic KdV equations, since the advantages of two kinds of well-posedness frame were carefully compared.

Because the proof of the averaging principle of \eqref{moskdv} has close correspondence with  well-posedness results of stochastic KdV equations, let us  revisit the relative well-posedness results in the following.

The deterministic KdV equation, which was first proposed in \cite{KdV}, is a well-known dispersive model to describe the motion of long, unidirectional water waves on a shallow channel. The Cauchy problem of the deterministic KdV equation in $H^s(\mathbb{R}),\ s\in \mathbb{R} $ is always an important problem in the fields of dispersive equations, harmonic analyses and integrable systems.  By an energy method, \cite{kdv_energy} proved the case $s>3/2$. Basing on smoothing effects of $U(t)$, \cite{dkdv} proved  the global well-posedness for  $s>3/4$.  Using the Bourgain space, \cite{bourgain} and \cite{kdv_m_3_4} proved cases $s=0$ and $s>-3/4$. By introducing a new workspace, \cite{guo2009global} proved the  well-posedness for the case $s=-3/4$. By the existence of commuting flows for integrable PDE, \cite{kdv_m_1} proved the global well-posedness in $H^{-1}(\mathbb{R})$. 

The stochastic KdV equations, which  describes the dispersive model with stochastic forces or in stochastic mediums (see \cite{forced_kdv})  is  also  widely investigated. For the additive stochastic KdV equations, \cite{KdVH1} and \cite{de1999white} proved the well-posedness in $L^2\left(\Omega;C([0,T];H^1(\mathbb{R}))\right)$ and $     L^2\left(\Omega;C([0,T];L^2(\mathbb{R}))\right),$ respectively.  Here, we should emphasize that in \cite{de1999white}, it introduced a kind of Bourgain space with time regularity index less than $1/2$ and low frequency restriction, which is the foundation of \cite{de1999white},  \cite{kdvmuti} and \cite{kdv_sharper}.

For multiplicative cases, by using a modified solution, \cite{kdvmuti} got the global well-posedness in $L^2\left(\Omega;C([0,T];L^2(\mathbb{R}))\right)$. In \cite{kdv_sharper}, by improving the bilinear estimate proposed in \cite{de1999white}, they released the restrictions of noises kernels. What's more, they also  overcame regularity problems in using conservation laws in $H^s(\mathbb{R}), s\in\mathbb{N}^+$, when prove the global well-posedness in $H^s(\mathbb{R}), s\in\mathbb{N}^+$, by an approximation approach.  In \cite{kdvas}, different from \cite{de1999white}, \cite{kdvmuti} and \cite{kdv_sharper}, they deal with the bilinear estimate with the workspace inspired by \cite{guo2009global}. As a result, in the local well-posedness proof of \cite{kdvas}, it no longer needs the homogeneous low frequency restriction of the workspace and the integrable conditions of noises kernels. When they prove the global well-posedness, the solution is extended only by the $L_x^2$ conservation law of the KdV equation. Thus, they can prove the global well-posedness in $ H^s(\mathbb{R}), s\geq 0$.

Now, in order to state the main theorems and prove the well-posedness of \eqref{moskdv} in $H^s(\mathbb{R})\times H^s(\mathbb{R})$, we will introduce the following notations and definitions.

For two positive constants $a,b$, we use $a\lesssim b$ ($ a\gtrsim b$) to represent that there exists a positive constant $C$ such that $ a\leq C b$ $(a>Cb)$. We also use $a\approx b$ to represent $a\lesssim b$ and $b\lesssim a$. What's more, let us set  $\left<C\right>=1+|C|$.

For any $T>0,\ p\in[1,\infty]$, 
$L^p_{T}:=\left\{f:\left(\int_{0}^{T}|f(t)|^p~dt\right)^{1/p}<\infty\right\} $ is the set of time integrable functions up to $T$. We also use $L_t^p$, $L_x^p$ to denote the set of time and space integrable functions on $\mathbb{R}$, respectively.  $C_TH^s_x, C_tH^s_x$ are notations of spaces $ C([0,T];H^s_x), C([0,\infty);H^s_x)$.

For any $g\in \mathscr{S}'(\mathbb{R}^2)$, we denote $\hat{g}(\tau,\xi), \mathscr{F}_{x}g(t,\xi)$  the space-time and the space Fourier transforms of $g$, respectively.
Correspondingly, $\check{g}(t,x), \mathscr{F}^{-1}_{\xi}g(t,x)$ are the inverse Fourier transforms of $g$.

For any $n\in\mathbb{N}^+$, $ f\in  \mathscr{S}'(\mathbb{R})$, we will use notations 
$$
P_nf(x):=\mathcal{F}^{-1}_\xi \left(\theta(\xi/n)\mathcal{F}_xf(\xi)\right),\ P_{> n}f(x):=\mathcal{F}^{-1}_\xi \left((1-\theta(\xi/n))\mathcal{F}_xf(\xi)\right)
$$ to represent the smooth low-frequency and the high-frequency cut-off operators, respectively. Here, $\theta\in C_0^\infty (\mathbb{R})$ and $\mathrm{supp}~\theta\subset [-2,2]$, $\theta|_{[-1,1]} = 1$.

Let $H^s(\mathbb{R}),\ s\in\mathbb{R}$ (simply denoted by $H_x^s$) be the normal Sobolev spaces with norms 
$$
\|f\|_{H_x^s}=\left(\int_{\mathbb{R}} |J^sf(x)|^2 dx \right)^{1/2}:=\left(\int_{\mathbb{R}_{\xi}} \left<\xi^2\right>^{s}\left|\mathscr{F}_{x}f(\xi)\right|^2 d\xi\right)^{1/2}.
$$
The homogeneous Sobolev norms are
$$
\|f\|_{\dot{H}_x^{q}}=\left(\int_{\mathbb{R}} |D^{q}f(x)|^2 dx \right)^{1/2}:=\left(\int_{\mathbb{R}_{\xi}} |\xi|^{2q}|\mathscr{F}_{x}f(\xi)|^2 d\xi\right)^{1/2}.
$$

For any $u\in H_x^s,s\geq 0$, we will use the following equality gotten by the Fubini Theorem and  the Leibniz-type estimates in \cite{benea2016multiple}
\begin{equation}\label{auxiliary_equality}
\begin{aligned}
\sum_{j=0}^{\infty}\|u\Phi_i e_j\|^2_{H_x^s}=&\int_{\mathbb
{R}_\xi}\left<\xi\right>^{2s}\sum_{j=0}^\infty\left|  \mathscr{F}_x[u\Phi_ie_j](\xi)\right|^2d\xi\\
=&\int_{\mathbb{R}_\xi}\left<\xi\right>^{2s}\int_{\mathbb{R}}\left|  \mathscr{F}_x[uk_i](\xi,y)\right|^2dy d\xi\\
=&\|J_x^s(uk_i)\|^2_{L_{x,y}^2}\\
\lesssim&\|u\|^2_{H_x^s}\|k_i\|^2_{L_x^\infty L_y^2}+\|u\|^2_{L_x^2}\|J^sk_i\|^2_{L_x^\infty L_y^2}\\
\lesssim&\|u\|^2_{H_x^s}\|J_x^sk_i\|^2_{L_x^\infty L_y^2},\ i=1,2.
\end{aligned}
\end{equation}

According to the Minkowski inequality, we have 
$\left\|J_x^sP_nk(x,y)\right\|_{L_x^\infty L_y^2}  \lesssim C_s\langle  n\rangle^s\|k(x,y) \|_{L_x^\infty L_y^2} $, for any $s\geq 0$.

For any $f\in \mathscr{S}'(\mathbb{R}^2)$, $b,s\in\mathbb{R}$, let us denote the Bourgain norms  and the Bourgain norm with low-frequency homogeneous restriction by
$$
\|f\|_{Y_{b,s}}=\left(\int_{\mathbb{R}^2}\left<\xi\right>^{2s} \left<\tau-\xi^3\right>^{2b}|\hat{f} (\xi,\tau)|^2d\xi d\tau\right)^{1/2}.
$$
and 
$$
\|f\|_{Y_{b,s,-q}}=\left(\int_{\mathbb{R}^2}|\xi|^{-2q}\left<\xi\right>^{2s} \left<\tau-\xi^3\right>^{2b}|\hat{f} (\xi,\tau)|^2d\xi d\tau\right)^{1/2},
$$
respectively.

Correspondingly, for any $T>0$, the restricted norms are denoted  by
$$
\|f\|_{Y^T_{b,s}}:=\inf\{\|\tilde{f}\|_{Y_{b,s}}:\tilde{f}=f,\ \forall t\in[0,T] \}.
$$
and 
$$\|f\|_{Y^T_{b,s,-q}}:=\inf\{\|\tilde{f}\|_{Y_{b,s,-q}}:\tilde{f}=f,\ \forall t\in[0,T] \},
$$
respectively.

When we establish an averaging principle converging in the sense of probability, we will use  $X^T_{b,s}$, which is defined by 
$$
\|f\|_{X^T_{b,s}}:=\|f\|_{Y^T_{b,s}}+ \|P_1f\|_{L_x^2L_T^\infty},
$$
as the workspace of the stochastic KdV type equation.

What's more, when we establish the general type averaging principle converging in the sense of $L^2_\omega$, we will use workspace $Z^{T}_{b,0,q}=Y^T_{b,0,-q}\cap Y^T_{b,0}$ equipped with the norm 
$$
\|f\|_{Z_{b,0,q}^T}=\|f\|_{Y_{b,0}^T}+\|f\|_{ Y_{b,0,-q}^T},
$$
as the workspace of the stochastic KdV type equation.

For $v_\delta$, although we need to deal with the multilinear estimate of the coupling term in KdV type equation, the $\|\cdot\|_{L_T^{p}H_x^s}$ norm is enough. What's more, because $S(t)$ is not a unitary group of $H_x^s$, there will be problems when we estimate the stochastic convolution terms, if we use other norms for dispersive equations (see for example \cite{kdvmuti}).

The main theorems of this paper are listed in the following.
The first theorem provides the almost surely global well-posedness of \eqref{moskdv}.

\begin{thm}\label{thm:global_well-posedness}
	Let $T_0>0$ and $s\geq 0$ be two arbitrary fixed constants. Under the conditions  
	$$u_0\in H_x^s,\ \|J_x^s k_1\|_{L_x^\infty L_y^2}<\infty,\ v_0\in H_x^s,\ \|J_x^sk_2\|_{L_x^\infty L_y^2}<\infty,
	$$ 
	we have the unique solution $ (u_\delta,v_\delta)$ of \eqref{moskdv} satisfies $u_\delta\in  C([0,T_0];H_x^s)\cap X^{T_0}_{b,s},\ \forall b\in(5/12,1/2)$ almost surely and $v_\delta\in L^{2m}(\Omega;L^{2m}_{T_0} H_x^s),\ \forall m\in\mathbb{N}^+$, for any  $\delta\in(0,1]$. 
\end{thm}

Different from the classical averaging principle of the SPDE, the following theorem establishes convergence not in the sense of $L^p_\omega$ but in the sense of probability, as $\delta \downarrow 0$. Thus, we can call it a weak type averaging principle. 

\begin{thm}\label{thm:average_principle}
Under the conditions of Theorem \ref{thm:global_well-posedness}, let us further assume that $\|J_x^sk_2\|_{L_x^\infty L_y^2}$ is less than the constant $ C^{\dagger}$ mentioned in Lemma \ref{lem:well-posedness_dissipative}. Then, for any $\zeta\in(0,1)$, $\alpha>0$,  there exists a  $\delta_0\in(0,1)$  such that 
$$
\mathbb{P}\left(\|u-u_{\delta}\|_{C([0,T_0];H_x^s)}>\alpha\right)\leq \zeta,\ \forall \delta\in(0,\delta_0),
$$
which means $\{u_\delta\}$ converges to $u$ in probability. Here, $u$ is the solution of the following averaging equation
\begin{equation}\label{average_eq}
	\left\{
	\begin{aligned}
		du&=-\left(\partial_{xxx} u+\frac{1}{2}\partial_x(u^2)\right)dt+u\Phi_1dW_t^1;\\
		u&=u_0.
	\end{aligned}\right.
\end{equation}
\end{thm}
\begin{rem}\label{rem:why_no_rate_classical}
	Different from Theorem \ref{thm:result_with_classical_workspace}, in Remark \ref{rem:convergence_rate_probability}, we have explain why we can not get the convergence rate of Theorem \ref{thm:average_principle}.
\end{rem}

The next theorem prove the averaging principle in the  classical sense, a convergence rate is also gotten.

\begin{thm}\label{thm:result_with_classical_workspace}
	Let $k_1(x,y)=k_1(x-y)$(Namely, $\Phi_1$ is a convolution type operator). For any fixed $T_0>0$,  $\forall  q\in (0,2/19)$, $ \forall \delta\in(0,1)$, if we have $k_1\in L_x^2\cap L_x^{2/(1+4q)}$ and  $\|k_2\|_{L_x^\infty L_y^2} $ less than the constant $C^\dagger(a,c)$ in Proposition \ref{prop:well-posedness_with_moments}, then we have 
	\begin{equation*}\label{classical_convergence}
	\lim_{\delta\downarrow0}\mathbb{E}\|u-u_\delta\|^2_{C([0,T_0];L_x^2)}=0,
	\end{equation*}
	where $u$ is the solution of \eqref{average_eq}.
	
	As for the convergence rate , according to \eqref{converging_speed}, for any $l>1$, there exists a common constant $C>0$ and a $\delta_{0}\lesssim 1/l \cdot C^{-\left(l^{1/(1-2a)}\vee l^{1/(2-2(a+c))}\right)}$ such that
	$$
	\mathbb{E}\|u-u_\delta\|^2_{C([0,T_0];L_x^2)}\lesssim 1/l,\ \forall \delta\in(0,\delta_0).
	$$
\end{thm}

When we prove averaging principles, the differences between two kinds of well-posedness frameworks are compared in the following remark.   

\begin{rem}\label{rem:explain_different_workspace_and_result}
In the proofs of Theorem \ref{thm:average_principle} and Theorem \ref{thm:result_with_classical_workspace}, we can see the different advantages and disadvantages of the different well-posedness frameworks for stochastic KdV equation. 

In Theorem \ref{thm:result_with_classical_workspace}, we use the workspace proposed in \cite{kdvmuti} and improved in \cite{kdv_sharper}. For this case, we can get an averaging principle in the $L^2_\omega(C([0,T_0];L_x^2))$ sense with convergence rate. However, the noise kernel need  the integrable condition $k_1\in L_x^{2/(1+4q)}$, which is not parallel to the deterministic case. More important, in the proof of Proposition \ref{prop:well-posedness_with_moments}, when we  extend the solution in $H_x^s,\ s \in \mathbb{N}$, we need to use the  $H_x^s,\ s \in \mathbb{N}$ conservation structure of the coupled  KdV type equation. Unfortunately, it is not easy to keep  $H_x^s,\ s \in \mathbb{N}$ conservation laws for the coupled equation. Thus, for \eqref{moskdv}, only the  $L^2_\omega(C([0,T_0];L_x^2))$  convergence can be gotten owing to the conservation structure.

In Theorem \ref{thm:average_principle}, we use the  workspace proposed in \cite{kdvas} to deal with the averaging principle problem. On the positive side, we only need the regularity conditions parallel to deterministic cases. More important, in the proof of Theorem \ref{thm:global_well-posedness}, the $L_x^2$ conservation law is enough to help us extend the solution in $H_x^s,\ s\geq 0$. Thus, when we study other dynamical problems in $H_x^s,\ s\geq 0$ of stochastic KdV equations, we can only  require the $L_x^2$ conservation structure of \eqref{moskdv}, which will provide considerable flexibility in selecting the model structure. On the negative side, referring to the explanation in the introduction of \cite{kdvas},  only the almost surely well-posedness and the weak type averaging principle with no convergence rate can be gotten.  

These  observations may also be important, when we treat  dynamical properties of general stochastic dispersive equations. 
\end{rem}

This paper is organized in the following manner:
In Section \ref{sec:esti}, we develop some necessary estimates for proving the local well-posedness of \eqref{moskdv}.   Next, in Section \ref{sec:global_well-posedness}, we get the local well-posedness of \eqref{moskdv} by a fixed point argument and then extend the solution in $H_x^s,\ s\geq 0$ by a $L_x^2$ prior estimate. The  key point of this section is that we can get a common moment estimate  for the solution up to a common small time. This moment estimate is independent of the initial value, as long as we have a prior estimate of $L_\omega^2L_x^2$ norm of initial values. In Section \ref{sec:averaging_principle}, we prove an averaging principle in the sense of probability, by a uniform tail probability estimate, the property of the averaging equation and the convergence of $u_{\delta,R}$. Finally, in Section \ref{sec:averaging_principle_with_moment}, by convergences before a kind of stopping times and a uniform  tail moment estimates, we prove an averaging principle with convergence rate in the $L_\omega^2$ sense.

\section{Fundamental Inequalities}\label{sec:esti}
We will first introduce some fundamental inequalities, which are necessary in dealing with  linear terms, nonlinear terms and stochastic convolutions. 

The following lemma about Bourgain-type norms is mentioned in \cite{erdogan_tzirakis_2016}.
\begin{lemma}(\cite{erdogan_tzirakis_2016})\label{lem:norm_inequalities}
The following inequality 
$$
\|f\|_{Y^T_{b',s}}\lesssim T^{b-b'}\|f\|_{Y^T_{b,s}},
$$
holds
for any $T>0$, $s\in\mathbb{R}$ and $-1/2<b'<b<1/2$.
\end{lemma}

The following restricted estimate is necessary in the proof of local well-posedness.
\begin{lemma}(\cite{erdogan_tzirakis_2016})\label{lem:restricted_estimate}
	The following inequality
	\begin{equation}\label{restricted_kdv}
		\left\|\int_{0}^{t}U(t-r)f dr\right\|_{Y^T_{b,s}}\lesssim \| f\|_{Y^T_{b',s}}
	\end{equation}
	holds for any $T\in[0,1]$, $s\in\mathbb{R}$, $b\in(1/2,1) $ and $ b=b'+1$.
\end{lemma}

\begin{coro}\label{coro:restricted_estimate}
	The restricted estimates 
	$$
	\left\|\int_{0}^{t} U(t-r) fdr\right\|_{Y^T_{c,s}}\lesssim T^{1-a-c}\|f\|_{Y^T_{-a,s}},\ 
	$$
	$$
	\left\|\int_{0}^{t} U(t-r) fdr\right\|_{Y^T_{b,s}}\lesssim T^{1-a-b-\varepsilon}\|f\|_{Y^T_{-a,s}} 
	$$
	hold for any $T\in[0,1]$, $s\in\mathbb{R}$, $a,b\in (0,1/2)$, $c\in(1/2,1)$, $\varepsilon\in(1/2-b,1-a-b)$ and satisfying $a+c<1$.
\end{coro}
\begin{proof}
The first formula is obvious. For the second one, we have 
$$
\left\|\int_{0}^{t} U(t-r) fdr\right\|_{Y^T_{b,s}}\lesssim \left\|\int_{0}^{t} U(t-r) fdr\right\|_{Y^T_{b+\varepsilon,s}} \lesssim T^{1-a-b-\varepsilon}\|f\|_{Y^T_{-a,s}},
$$
by  Lemma \ref{lem:norm_inequalities} and Lemma \ref{lem:restricted_estimate}.
Here, we need $b+\varepsilon>1/2$ and $b+\varepsilon+a<1$ to fit restrictions in those lemmas.
\end{proof}

In the rest of this paper, when we choose a small interval $[0,T]$, we assume $T\in(0,1]$ by default.

The bilinear estimates of the KdV type equation have been gotten in \cite{kdvas}. 
\begin{lemma}(\cite{kdvas})\label{lem:kdv_bilinear_esti}
	Let $1/4<a,b<1/2$, $a+2b\geq 5/4$, $s\geq 0$ and $0<T<1$. We have
	\begin{equation}\label{bilinearclear}
		\begin{aligned}
			\|(uv)_x\|_{Y^T_{-a,s}}\lesssim& \Big(\|P_1 u\|_{L_T^\infty L_x^2}\|P_1v\|_{L_{T}^2L_x^\infty}+ \|P_1 u\|_{L_x^2L_T^\infty}\|P_{>1}v\|_{Y^T_{b,s}} \\
			&+ \|P_{>1}u\|_{Y^T_{b,s}}\|P_1 v\|_{L_x^2L_T^\infty}+ \|P_{>1}u\|_{Y^T_{b,s}}\|P_{>1}v\|_{Y^T_{b,0}}\\
			&+\|P_{>1}u\|_{Y^T_{b,0}}\|P_{>1}v\|_{Y^T_{b,s}}
			\Big),
		\end{aligned}
	\end{equation}
	which means 
	\begin{equation*}\label{bilinearest}
		\|(uv)_x\|_{Y^T_{-a,s}}\lesssim \|u\|_{X^T_{b,s}}\|v\|_{X^T_{b,0}}+\|u\|_{X^T_{b,0}}\|v\|_{X^T_{b,s}}.
	\end{equation*}
\end{lemma}

\begin{lemma}(\cite{kdvas})\label{lem:lowfrebi}
	For any $0<T<1$, we have
	\begin{equation*}
		\begin{aligned}
			&\left\|P_1\int_{0}^{t} U(t-t')(u v)_x~d{t'}\right\|_{L_x^2 L^\infty_T}\\
			\lesssim &T^{1/2}\left(\|P_1u\|_{L_x^2L_T^\infty} +\|P_{>1}u\|_{Y_{1/3,0}^T}\right)\left(\|P_1v\|_{L_x^2L_T^\infty} +\|P_{>1}v\|_{Y_{1/3,0}^T}\right).
		\end{aligned}
	\end{equation*}
 
\end{lemma}

In the proof of the trilinear estimate of the KdV type equation, we will use the following lemma.

\begin{lemma}\label{lem:interpolation_bourgain_norm_esti}
	The following estimate 
	\begin{equation}\label{interpolation_bourgain_norm_esti}
		\|f\|_{L_T^4L_x^\infty}\lesssim \|f\|_{Y^T_{b,0}}
	\end{equation}
holds for any $T>0$ and $b\in(5/12,1/2)$.
\end{lemma}
\begin{proof}
According to the proof of Strichartz estimates for KdV equation in \cite{dkdv} and the transfer principle in \cite{Tao2006}, we have  
\begin{equation}\label{kdv_bourgain_esti_1}
\left\|J_x^{3/p-1/2}P_{>1}f\right\|_{L_t^pL_x^\infty}\approx\left\|D_x^{3/p-1/2}P_{>1}f\right\|_{L_t^pL_x^\infty}\lesssim \|f\|_{Y_{b,0}},
\end{equation}
for any $p\in[4,6]$ and $ b>1/2$.

By the Bernstein inequality, the maximal functional estimate in \cite{dkdv} and the transfer principle, we have 
\begin{equation}\label{kdv_bourgain_esti_2}
	\begin{aligned}
	\left\|J_x^sP_{1}f\right\|_{L_t^pL_x^\infty}\approx&\|P_1f\|_{L_t^pL_x^\infty}
		\lesssim \|P_1f\|_{L_t^pL_x^2}\\
		\lesssim&\|P_1f\|^{1-2/p}_{L_t^\infty L_x^2}\|P_1f\|^{2/p}_{L_{t,x}^2}\\
		\lesssim&\|P_1f\|_{Y_{b,0}},\ \forall s\in[0,1/4],b>1/2, 1\leq p\leq \infty.
	\end{aligned}
\end{equation} 

Thus, combining \eqref{kdv_bourgain_esti_1} with \eqref{kdv_bourgain_esti_2}, we have 
\begin{equation}\label{kdv_bourgain_esti}
	\begin{aligned}
\left\|f\right\|_{L_t^pL_x^\infty}\lesssim \left\|P_1f\right\|_{L_t^pL_x^\infty}+\left\|P_{>1}f\right\|_{L_t^pL_x^\infty}\lesssim  \|f\|_{Y_{b,1/2-3/p}}, \ \forall p\in[4,6],\ b>1/2.
	\end{aligned}
\end{equation}

What's more, by the Sobolev embedding Theorem,
we have 
$$
\|f\|_{L^p_tL_x^\infty}\lesssim \|f\|_{H_t^bH_x^s}, \  b\geq 1/2-1/p,\ s>1/2,
$$
for any $ p\geq2$, which means 
$$
\|f\|_{L^p_tL_x^\infty}\lesssim \|f\|_{Y_{b,s}}, \  b\geq 1/2-1/p,\ s>1/2,
$$
because of $\|f\|_{H_t^bH_x^s}=\|U(-t)f\|_{H_t^bH_x^s}=\|f\|_{Y_{b,s}} $.

Thus, we have
$$
\|f\|_{L_t^4L_x^\infty}\lesssim \|f\|_{Y_{b_1,s_1}},\ \|f\|_{L_t^4L_x^\infty}\lesssim \|f\|_{Y_{b_2,s_2}},
$$
for $s_1=-1/4$ and $\forall b_1>1/2, \forall s_2>1/2, \forall b_2>1/4$.

To get \eqref{kdv_bourgain_esti}, we use the interpolation space 
$$Y_{\theta b_1+(1-\theta)b_2,\ \theta s_1+(1-\theta)s_2},\  \theta\in[0,1].$$
Let $ \theta s_1+(1-\theta)s_2= 0$. By the Stein interpolation  in \cite{stein_functional}, we have 
$$
\|f\|_{L_t^4L_x^\infty}\lesssim \|f\|_{Y_{b,0}},\ \forall b\in(5/12,1/2).
$$

According to the definition of $\|\cdot\|_{Y^T_{b,0}}$, for any $ \varepsilon>0$, there exists a $g\in Y_{b,0}$ such that $g=f,\ t\in[0,T] $ and $ \|g\|_{Y_{b,0}}\leq \|f\|_{Y^T_{b,0}}+\varepsilon$. Therefore, we have
\begin{equation}\label{auxiliary}
	\|f\|_{L_T^4L_x^\infty}\leq \|g\|_{L_t^4L_x^\infty}\leq \|g\|_{Y_{b,0}}\leq \|f\|_{Y^T_{b,0}}+\varepsilon,
\end{equation}
which finishes the proof.
\end{proof}

Now, we prove the trilinear estimate of the nonlinear term in the KdV type equation.

\begin{lemma}\label{lem:tri_esti_kdv_1}
	The following inequality
	\begin{equation}\label{tri_esti_kdv_1}
		\begin{aligned}
	\left\|fgh\right\|_{Y^T_{-a,s}}\lesssim& \|f\|_{Y^T_{b,s}}\|g\|_{L_T^4L^2_x}\|h\|_{L_T^4L^2_x}+\|f\|_{Y^T_{b,0}}\|g\|_{L_T^4H^s_x}\|h\|_{L_T^4L^2_x}\\
	&+ \|f\|_{Y^T_{b,0}}\|g\|_{L_T^4L^2_x}\|h\|_{L_T^4H^s_x},
		\end{aligned}
	\end{equation}
holds for any $T>0$, $s\geq 0$, $a,b\in(5/12,1/2)$.
\end{lemma}
\begin{proof}
Just like \eqref{auxiliary}, we  verify
$$
\begin{aligned}
\left\|fgh\right\|_{Y_{-a,s}}\lesssim&\|f\|_{Y_{b,s}}\|g\|_{L_t^4L^2_x}\|h\|_{L_t^4L^2_x} +\|f\|_{Y_{b,0}}\|g\|_{L_t^4H_x^s}\|h\|_{L_t^4L^2_x}\\
&+\|f\|_{Y_{b,0}}\|g\|_{L_t^4L^2_x}\|h\|_{L_t^4H_x^s},\ \forall a,b\in(5/12,1/2).
\end{aligned}
$$

By duality, we only need to prove
$$
\begin{aligned}
(v,fgh)
\lesssim &\|v\|_{Y_{a,-s}}\big( \|f\|_{Y_{b,s}}\|g\|_{L_t^4L^2_x}\|h\|_{L_t^4L^2_x} +\|f\|_{Y_{b,0}}\|g\|_{L_t^4H_x^s}\|h\|_{L_t^4L^2_x}\\
&+\|f\|_{Y_{b,0}}\|g\|_{L_t^4L^2_x}\|h\|_{L_t^4H_x^s}\big),\ \forall a,b\in(5/12,1/2).
\end{aligned}
$$

According to the H\"{o}lder inequality and the  Leibniz-type estimates in \cite{benea2016multiple}, we have 
\begin{equation}\label{counter_example}
\begin{aligned}
	&\left(v,fgh\right)\\
	\lesssim&\|J_x^{-s}v\|_{L_t^4L^\infty_x}\|J_x^{s}(fgh)\|_{L_t^{4/3}L^1_x}\\
	\lesssim&\|v\|_{Y_{a,-s}}\big(\|J^s_xf\|_{L_t^4L_x^\infty}\|g\|_{L_t^4L_x^2}\|h\|_{L_t^4L_x^2}+\|f\|_{L_t^4L_x^\infty}\|g\|_{L_t^4H_x^s}\|h\|_{L_t^4L_x^2}\\
	&+\|f\|_{L_t^4L_x^\infty}\|g\|_{L_t^4L_x^2}\|h\|_{L_t^4H_x^s}\big).
\end{aligned}
\end{equation}

Thus, by \eqref{interpolation_bourgain_norm_esti}, we can prove \eqref{tri_esti_kdv_1}.
\end{proof}

\begin{lemma}\label{lem:tri_esti_kdv_2}
\begin{equation}\label{tri_esti_kdv_2}
	\left\|P_1\int_{0}^{t}U(t-r)(fgh)dr\right\|_{L_x^2L_T^\infty}\lesssim T^{1/4}\|f\|_{Y^T_{b,0}}\|g\|_{L_T^4L_x^2}\|h\|_{L_T^4L_x^2},
\end{equation}
for any $b\in(5/12,1/2)$.
\end{lemma}
\begin{proof}
By the Christ-Kiselev Lemma in \cite{molinet2004well} and the maximal functional estimate in \cite{dkdv}, we have 
$$
\begin{aligned}
	\left\|P_1\int_{0}^{t}U(t-r)(fgh)dr\right\|_{L_x^2L_T^\infty}\lesssim\left\| P_1 (fgh)\right\|_{L_T^1H_x^\sigma},
\end{aligned}
$$
for any $\sigma>3/4$.

Therefore, by the cut-off property in the frequency space and the Young inequality, we have 
\begin{equation*}
\begin{aligned}
	\left\| P_1 (fgh)\right\|_{L_T^1H_x^\sigma}\lesssim &\left\| P_1 (fgh)\right\|_{L_T^1L_x^2}\\
	\lesssim& \left\| fgh\right\|_{L_T^1L_x^1}\\
	\lesssim&T^{1/4}\|f\|_{L_T^4L_x^\infty}\|g\|_{L_T^4L_x^2}\|h\|_{L_T^4L_x^2}\\
	\lesssim&T^{1/4}\|f\|_{Y^T_{b,0}}\|g\|_{L_T^4L_x^2}\|h\|_{L_T^4L_x^2},
\end{aligned}
\end{equation*}
for any $b\in(5/12,1/2)$.

\end{proof}

In the rest part of this subsection, we will introduce the estimate of stochastic convolution terms in \eqref{moskdv}. 

\begin{lemma}\label{lem:SI_esti_dissipative}
Let $ m\in \mathbb{N}^+$, $s\geq0$ and $ \|J_x^sk_2\|_{L_x^\infty L_y^2} <+\infty$. If $ f\in L^{2m}_\omega L_T^{2m}H_x^s$, then the following estimate
\begin{equation}\label{SI_esti_dissipative}
\mathbb{E}\left\|\int_{0}^{t} S(t-r) (f \Phi_2 dW_r^2)\right\|_{L^{2m}_T H_x^s}^{2m}\lesssim T^{m}\|J_x^sk_2\|^{2m}_{L_x^\infty L_y^2}\mathbb{E}\| f\|_{L_T^{2m}L_x^2}^{2m},
\end{equation}
holds for any $T>0$.
\end{lemma}
\begin{proof}
By the H\"{o}lder inequality, the B\"{u}rkholder–Davis–Gundy(BDG) inequality and \eqref{auxiliary_equality}, we have  
$$
\begin{aligned}
&\mathbb{E}\left\|\int_{0}^tS(t-r)(f\Phi_2dW_r^2)\right\|^{2m}_{L_T^{2m}H_x^s}\\
\lesssim&\int_{0}^{T}\mathbb{E} \left\|\int_{0}^tS(t-r)(f\Phi_2dW_r^2)\right\|^{2m}_{H_x^s}dt\\
\lesssim&\int_{0}^{T}\mathbb{E}\sup_{t_0\in[0,t]} \left\|\int_{0}^{t_0}S(t-r)(f\Phi_2dW_r^2)\right\|^{2m}_{H_x^s}dt\\
\lesssim&T^{m}\|J_x^sk_2\|^{2m}_{L_x^\infty L_y^2}\mathbb{E}\| f\|_{L_T^{2m}H_x^s}^{2m}.
\end{aligned}
$$ 
\end{proof}

The following estimate of the stochastic convolution in KdV type equation comes from \cite{kdvas}.
\begin{lemma}(\cite{kdvas})\label{lem:estimate_of_SI_kdv}
Let $b\in(0,1/2),\ s\geq 0,\ T>0$. We have 
\begin{equation}\label{si_kdv}
	\begin{aligned}
		&\mathbb{E}\left\| \int_{0}^{t} U(t-r)(u\Phi_1dW_r^{1})\right\|^{2}_{X^{\sigma\wedge T}_{b,s}\cap C([0,\sigma\wedge T];H_x^s)}\\
		\lesssim&(T^{1-2b}+1)\|u\Phi_1e_j\|_{L^2_\omega L^2_{ \sigma \wedge T}l^2_j H^s_x}^2\\
	\lesssim &(T^{1-2b}+1) \|J_x^sk_1\|^{2}_{L_x^\infty L_y^2}\mathbb{E}\|u\|^{2}_{L_{\sigma\wedge T}^2H_x^s},
	\end{aligned}
\end{equation}
for any stopping time $\sigma\geq0$ almost surely.
\end{lemma}

\begin{rem}\label{rem:SI_esi_kdv}
From the proof of Lemma 3.12 in \cite{kdvas}, we   know that  only the $\|\cdot\|_{L_\omega^2X^T_{b,s}}$ norm rather than $ \|\cdot\|_{L_\omega^{2m}X^T_{b,s}}, m>1$ norm of the stochastic convolution in the KdV type equation can be estimated. Roughly speaking, this is because when we try to use the BDG inequality to estimate  $ L^{2m}_\omega L_t^2 L_x^\infty$-norm of the stochastic convolution terms, we need convert it to $  L_t^2 L^{2m}_\omega L_x^\infty $ and then back to $ L^{2m}_\omega L_t^2 L_x^\infty$. Thus, $m$ can only be $1$.
\end{rem}

The following two lemmas will be used in Section \ref{sec:averaging_principle_with_moment}, when we prove the averaging principle in the classical sense.
\begin{lemma}\label{lem:trilinear_estimate_q}
	The inequality 
	\begin{equation}\label{tri_esti_kdv_q}
		\begin{aligned}
			\left\|fgh\right\|_{Y^T_{-a,0,-q}}\lesssim  \|f\|_{Y^T_{b,0}}\| g\|_{L_T^4L_x^2} \|h \|_{L_T^4L_x^2},
		\end{aligned}
	\end{equation}
	holds for any $T>0$, $q\in [0,1/2]$, $ a\in(5/12-q/3,1/2) ,\ b\in (5/12,1/2)$.
\end{lemma}
\begin{proof}
	Just like Lemma \ref{lem:tri_esti_kdv_1}, it is sufficient to study 
	$$
	\begin{aligned}
		\left(v,fgh\right)=& (D^q(v),D^{-q}(fgh))\\
		\lesssim&\|D^qv\|_{L_T^4L_x^{1/q}}\|D^{-q}(fgh)\|_{L_T^{4/3}L_x^{1/(1-q)}}.
	\end{aligned}
	$$
	By the Sobolev embedding and \cite{strichartz_ariy}, we have 
	$$
	\begin{aligned}
		\left(v,fgh\right)\lesssim&\|D^qv\|_{L_T^4L_x^{1/q}}\|fgh\|_{L_T^{4/3}L_x^1}\\
		\lesssim&\|v\|_{Y^T_{a,0,q}} \|f\|_{Y^T_{b,0}}\| g\|_{L_T^4L_x^2} \|h \|_{L_T^4L_x^2}.
	\end{aligned}
	$$
	for any $ a\in(5/12-q/3,1/2) ,\ b\in (5/12,1/2)$.
\end{proof}

Based on the proofs of Proposition 2.5 and 2.7 in \cite{kdvmuti}, we have the following lemma.
\begin{lemma}\label{lem:stochastic_convolution_estimate_2}
	For any $m\in\mathbb{N}$, $T>0$, $b\in(0,1/2)$, $\varepsilon\in(0,1-2b)$ and $q>0$,  we have 
	\begin{equation}\label{stochastic_convolution_estimate_2.1}
		\begin{aligned}
			&\mathbb{E}\left\|\int_{0}^{t} U(t-r)\left(f\Phi_1 dW^{(1)}_r\right)\right\|^{2m}_{Z^T_{b,0,q}}\\
			\lesssim& T^{m(1-2b-\varepsilon)}(\|k_1\|^{2m}_{L_x^{2/(1+4q)}}+\|k_1\|_{L_x^2}^{2m})\mathbb{E}\|f\|^{2m}_{L_{T,x}^2}
		\end{aligned}
	\end{equation}
	and 
	\begin{equation}\label{stochastic_convolution_estimate_2.2}
		\begin{aligned}&\mathbb{E}\left\|\int_{0}^{t} U(t-r)\left(f\Phi_1 dW^{(1)}_r\right)\right\|^{2m}_{C([0,T];L_x^2\cap\dot{H}_x^{-q})}\\
			\lesssim& (\|k_1\|^{2m}_{L_x^{2/(1+4q)}}+\|k_1\|_{L_x^2}^{2m})\mathbb{E}\|f\|^{2m}_{L_{T,x}^2}.
		\end{aligned}
	\end{equation}
\end{lemma}

\section{The well-posedness of \eqref{moskdv}}\label{sec:global_well-posedness}

Firstly, let us consider the well-posedness of the stochastic dissipative equation.

\begin{lemma}\label{lem:well-posedness_dissipative}
Let $T_0>0$, $s\geq0$ be two settled constants. Suppose that $v_0\in H_x^s$ and  $\|J_x^sk_2\|_{L_x^\infty L_y^2}<\infty$. Then, $v_{\delta}$ is well-posed in $ L^{2m}(\Omega;L^{2m}_{T_0}H_x^s)$, $\forall m\in\mathbb{N}^+$. 

Moreover, for any $m\in\mathbb{N}^+$, there exists a constant $C^{\dagger}(m)>0$, such that if $\|J_x^sk_2\|_{L_x^\infty L_y^2}\leq C^{\dagger}(m)$, then 
\begin{equation}\label{decrease_delta}
	\|v_{\delta}\|^{2m}_{L_\omega^{2m}L_{T_0}^{2m}H_x^s}\lesssim C(\|v_0\|_{H_x^s},T_0)\cdot m^{-1} \delta.
\end{equation}
\end{lemma} 
\begin{proof}
	We set 
	$$
	\mathcal{T} v_{\delta}:=S(t)v_{0} -\delta^{-1}\int_{0}^{t} S(t-r) v_\delta dr+\delta^{-1/2}\int_{0}^{t} S(t-r) (v_\delta \Phi_2 dW_r^2).
	$$
	By Lemma \ref{lem:SI_esti_dissipative}, it is clear that 
	$$
	\begin{aligned}
	&\left\|\mathcal{T}  v_{\delta}\right\|_{L_\omega^{2m}L_T^{2m}H_x^s}\\
	\lesssim &C(T)\|v_0\|_{H_x^s} +(\delta^{-1}T+\delta^{-1/2}T^{1/2} \|J_x^sk_2\|_{L_x^\infty L_y^2})\left\| v_{\delta}\right\|_{L_\omega^{2m}L_T^{2m}H_x^s}
	\end{aligned}
	$$
	and
	$$
	\begin{aligned}
	&\left\|\mathcal{T}  v^{(1)}_{\delta}-\mathcal{T}v^{(2)}_{\delta}\right\|_{L_\omega^{2m}L_T^{2m}H_x^s}\\
	\lesssim& (\delta^{-1}T+\delta^{-1/2}T^{1/2} \|J_x^sk_2\|_{L_x^\infty L_y^2})\left\| v^{(1)}_{\delta}- v^{(2)}_{\delta}\right\|_{L_\omega^{2m}L_T^{2m}H_x^s},
	\end{aligned}
	$$
	for any $T\in(0,1)$, $ m\geq 1$. Thus, if we choose $T^*(\delta,\|J_x^sk_2\|_{L_x^\infty L_y^2})$ sufficiently small, then $\mathcal{T}$ is bounded and contractive in $L_\omega^{2m}L_{T^*}^{2m}H_x^s$.
		
	By the It\^{o} formula, we have
	$$
	\begin{aligned}
	&(J_x^sv_{\delta},J_x^sv_{\delta})^{m}\\
	\leq&(J_x^s v_0,J_x^s v_0)^m-2m\delta^{-1} \int_{0}^{t}(J_x^sv_{\delta}, J_x^sv_{\delta})^mdr\\
	&+2m\delta^{-1/2}\int_{0}^{t}(J_x^sv_{\delta}, J_x^sv_{\delta})^{m-1}\left(J_x^s v_{\delta}, J_x^s(v_{\delta}\Phi_2dW_r^{2})\right)\\
	&+m\delta^{-1}\sum_{i=0}^{\infty}\int_{0}^{t}(J_x^sv_{\delta}, J_x^sv_{\delta})^{m-1} \|J_x^s\left(v_{\delta}\Phi_2e_i\right)\|^2_{L_x^2} dr\\
	&+2m(m-1)\delta^{-1}\sum_{i=0}^{\infty}\int_{0}^{t}(J_x^sv_{\delta}, J_x^sv_{\delta})^{m-2}\left(J_x^s v_{\delta}, J_x^s(v_{\delta}\Phi_2e_i)\right)^2dr,
	\end{aligned}
	$$
	for any $ m\in\mathbb{N}^+,\ t\in[0,T^*]$.

	Thus, by the Gronwall inequality, we have 
	$$
	\mathbb{E}\|v_{\delta}(t)\|_{H_x^s}^{2m}\leq C(m,T_0,\delta,\|v_0\|_{H_x^s},\|J_x^sk_2\|_{L_x^\infty L_y^2}),\ \forall t\in[0,T^*].
	$$
	
	Therefore, by dividing $[0,T_0]$ into finite numbers of small intervals 
	$$[0,T^*],[T^*,2T^*],...,[T_0-T^*,T_0],$$
	and using fixed point arguments, we can prove $v_\delta$ is well-posed in $L_\omega^{2m}L_{T_0}^{2m}H_x^s$, for any $m\geq 1$, $m\in\mathbb{N}^+$. 

	Next, we prove \eqref{decrease_delta}.
	For any fixed $m\in\mathbb{N} $, since $v_\delta$ is well-posed in $L_\omega^{2m}L_{T_0}^{2m}H_x^s$,  the It\^{o} formula in $[0,T_0]$ shows
	$$
	\begin{aligned}
	&\mathbb{E}(J_x^sv_{\delta},J_x^sv_{\delta})^m \\
	\leq&(J_x^s v_0,J_x^s v_0)^m-2m\delta^{-1} \mathbb{E}\int_{0}^{t}(J_x^sv_{\delta}, J_x^sv_{\delta})^mdr \\
	&+Cm\delta^{-1} \|J_x^sk_2\|^2_{L_x^\infty L_y^2}\mathbb{E}\int_{0}^{t}(J_x^sv_{\delta}, J_x^sv_{\delta})^{m} dr\\
	&+2Cm(m-1)\|J_x^sk_2\|^2_{L_x^\infty L_y^2}\delta^{-1}\mathbb{E}\int_{0}^{t}(J_x^sv_{\delta}, J_x^sv_{\delta})^{m}dr ,
	\end{aligned}
	$$
	almost surely for any $t\in[0,T_0]$.
	
	Therefore, there exists a  constant $C^{\dagger}(m)$ only depending on $m$ such that 
	\begin{equation}\label{esti_gronwall_auxi}
	\begin{aligned}
	&\mathbb{E}(J_x^sv_{\delta},J_x^sv_{\delta})^m \\
	\leq&(J_x^s v_0,J_x^s v_0)^m-m\delta^{-1} \mathbb{E}\int_{0}^{t}(J_x^sv_{\delta}, J_x^sv_{\delta})^mdr ,
	\end{aligned}
	\end{equation}
	once $ \|J_x^sk_2\|_{L_x^\infty L_y^2}\lesssim C^{\dagger}(m)$.
	
	Hence, by the Gronwall inequality, it can be deduced from \eqref{esti_gronwall_auxi} that 
	\begin{equation}\label{L_x^2_pri_esti_2}
	\begin{aligned}
		\|v_{\delta}(t)\|^{2m}_{L_\omega^{2m}H_x^s}\leq C(\|v_0\|_{H_x^s})e^{-m\delta^{-1} t},\ \forall t\in[0,T_0],
	\end{aligned}
\end{equation} 
which finishes the proof of  this lemma. 

\end{proof}

\begin{coro}\label{coro:longtime_dissipative_equation}
Under the condition of Lemma \ref{lem:well-posedness_dissipative}, according to \eqref{L_x^2_pri_esti_2} and Theorem 4.3.9  in \cite{liuweispde}( or Remark 2.1 in \cite{apabstract}), for any $v_0\in H_x^s$, $v_\delta$ has a unique invariant measure $\delta_0^{H_x^s}$. Here, $\delta_0^{H_x^s}$ means the Dirac measure concentrated on $0$ in $ H_x^s$.
\end{coro}

Next, by a fixed point argument, we will prove the local well-posedness of the KdV type equation in $ H_x^s$. To deal with the nonlinear terms, we need to introduce following notations.

Let us fixed a $ s\geq 0$. For any $ R>0$, we set
$$B_R(u):= \left(\theta_R(\|P_1u\|_{L_x^2L^\infty_t})P_1u+\theta_R\left(\|u\|_{Y^t_{b,0}}\right)P_{>1}u(x,t)\right)^2$$ 
and
$$v_{\delta,R}:=\theta_R(\|v_\delta\|_{L_t^4H_x^s})v_{\delta},\ \theta_R(\cdot):=\theta(\cdot^2/R^2).$$

Here,  $\theta(\cdot )$ is a smooth cut-off function supported on $[-2,2]$ and equal to $1$ in $[-1,1]$.

Let us introduce the following localized equation 
\begin{equation}\label{truncated_eq}
	\begin{aligned}
	u_{\delta,R}(t)=&U(t)u_0-1/2\int_{0}^{t} U(t-r)(B_R(u_{\delta,R}))_xdr\\
	&-\int_{0}^{t} U(t-r)(u_{\delta,R}v_{\delta,R}^2)dr\\
	&+\int_{0}^{t}U(t-r) (u_{\delta,R}\Phi_1 dW_{r}).
	\end{aligned}
\end{equation}

For any $\forall R,R_1>0$,  we also define the following stopping times
$$
\sigma^s_{\delta,R}(\omega):=\inf\{t\geq 0:\|v_{\delta}\|_{L_t^4H_x^s}\geq R\} \wedge T_0,$$ $$ \sigma_{\delta,R}(\omega):=\inf\{t\geq 0:\|v_{\delta}\|_{L_t^4L_x^2}\geq R\} \wedge T_0
$$
$$
\tau_{\delta,R}(\omega) := \inf\{t\geq 0:\|u_{\delta,R}\|_{X^t_{b,0}}\geq R \} \wedge T_0, 
$$
$$
\tau_{\delta,R,s}(\omega) := \tau_{\delta,R}(\omega) \wedge \sigma^s_{\delta,R},\ 
\tau_{\delta,R,R_1,s}(\omega) := \tau_{\delta,R}(\omega)\wedge \sigma^s_{\delta,R_1},
$$
almost surely.

The following lemma is the Lemma 3.2 in \cite{kdvas}, which will be used in the proof of local well-posedness.

\begin{lemma}(\cite{kdvas})\label{lem:truncated_equation}
	Let $b\in[0,1/2)$. For any $R>0$, $T>0$, $s\geq 0$, $f, g\in Y_{b,s}^T$, we have
	$$ 
	\left\|\theta\left(\|f\|^2_{Y^t_{b,0}}/R^2\right)f(x,t)\right\|_{Y_{b,0}^T}\lesssim R,
	$$
	$$
	\left\|\theta\left(\|f\|^2_{Y^t_{b,0}}/R^2\right)f(x,t)\right\|_{Y_{b,s}^T}\lesssim \|f\|_{Y_{b,s}^T},
	$$
	$$
	\left\|\theta\left(\|f\|^2_{Y^t_{b,0}}/R^2\right)f(x,t)-\theta\left(\|g\|^2_{Y^t_{b,0}}/R^2\right)g(x,t)\right\|_{Y^T_{b,0}}\lesssim \|f-g\|_{Y^T_{b,0}}.
	$$
\end{lemma}

In the following lemma, we will prove the local well-posedness of the KdV type equation.
\begin{lemma}\label{lem:local_well-posedness_kdv}
Let $T_0>0$, $s\geq0$ be two settled constants. Suppose that $\delta\in(0,1)$, $(u_0,v_0)\in H_x^s\times H_x^s$,  $\|J_x^sk_1\|_{L_x^\infty L_y^2}<\infty$ and $\|v_{\delta}\|_{L_{T_0}^4H_x^s}<\infty,\  \text{a.s.} \mathbb{P} $. Then, we have $u_{\delta,R}$ is  well-posedness in $L^2\left(\Omega; C([0,T_0];H_x^s)\cap X^{T_0}_{b,s}\right)$ for any $b\in(5/12,1/2)$. What's more, we have 
$$
\mathbb{E}\|u_{\delta,R}\|^2_{C([0,{T_0}];H_x^s)\cap X^{T_0}_{b,s}} \lesssim C(\|u_0\|_{H_x^s},R,T_0,\|J_x^sk_1\|_{L_x^\infty L_y^2}).
$$
\end{lemma}
\begin{proof}
	We first prove $u_{\delta,R}\in L^2\left(\Omega;C([0,T_0];L_x^2)\cap X_{b,0}^{T_0}\right)$. Let us set 
	$$
	\begin{aligned}
	\mathcal{T}u_{\delta,R}^{(i)}:=&U(t)u_0-\int_{0}^{t} U(t-r) \left(\frac{1}{2}\partial_{x}\left(B_R\left(u_{\delta,R}^{(i)}\right)\right)^2 +u^{(i)}_{\delta,R}v^2_{\delta,R}\right)dr\\
	&+\int_{0}^t U(t-r) \left(u_{\delta,R}^{(i)}\Phi_1 dW_r^1\right),\ i=1,2.
	\end{aligned}
    $$
    
	Thus, by Corollary \ref{coro:restricted_estimate}, Lemma \ref{lem:lowfrebi}, Lemma \ref{lem:kdv_bilinear_esti} and the proof of Lemma \ref{lem:tri_esti_kdv_2}, we have
	$$
	\begin{aligned}
	&\left\|\mathcal{T}u^{(1)}_{\delta,R}-\mathcal{T}u^{(2)}_{\delta,R}\right\|_{ C([0,T];L_x^2)\cap X^T_{b,0}}\\
	\lesssim& R\left(T^{1-(a+c)}+T^{1/4-a/2}+T^{1/2}\right)\left\|u^{(1)}_{\delta,R}-u^{(2)}_{\delta,R}\right\|_{X^T_{b,0}}\\
	&+R^2 (T^{1-(a+c)}+T^{1/4-a/2}+T^{1/4})\left\|u^{(1)}_{\delta,R}-u^{(2)}_{\delta,R}\right\|_{Y^T_{b,0}}\\
	&+\left\| \int_{0}^t U(t-r) \left(\left(u^{(1)}_{\delta,R}-u^{(2)}_{\delta,R}\right)\Phi_1 dW_r^1\right)\right\|_{C([0,T];L_x^2)\cap X^T_{b,0}},
	\end{aligned}
	$$
	for any $a,b\in(5/12,1/2)$, $c>1/2$, $a+c<1$, $T>0$.
	Here, to use Corollary \ref{coro:restricted_estimate}, we have chosen $\varepsilon=3/4-a/2-b$ for simplicity.  
	
	Hence, if we choose $T^*=T^*\left(R,\|k_1\|_{L_x^\infty L_y^2},a,b,c\right)>0$ sufficiently small, then by Lemma \ref{lem:estimate_of_SI_kdv}, we have 
	$$
	\begin{aligned}
		\left\|\mathcal{T}u^{(1)}_{\delta,R}-\mathcal{T}u^{(2)}_{\delta,R}\right\|_{L^2_{\omega}X^{T^*}_{b,0}}
		\lesssim& \frac{1}{2}\left\|u^{(1)}_{\delta,R}-u^{(2)}_{\delta,R}\right\|_{L_\omega^2X^{T^*}_{b,0}},\ \forall b\in(5/12,1/2).
	\end{aligned}
	$$
	
	Since it is easy to prove $\mathcal{T}$ is bounded in $ L_\omega^2\left( C([0,T^*];L_x^2)\cap X^{T^*}_{b,0}\right)$, we can prove $u_{\delta,R}$ is well-posed in $ L_\omega^2\left(C([0,T^*];L_x^2)\cap X^{T^*}_{b,0}\right)$ by a fixed point argument. Therefore, by dividing $[0,T_0]$ into finite numbers of small intervals shorter than $T^*$, we can prove $u_{\delta,R} $ is well-posed in $ L_\omega^2\left( C([0,T_0];L_x^2)\cap X^{T_0}_{b,0}\right)$ and 
	$$
    \|u_{\delta,R}\|^2_{L_\omega^2\left(C([0,T_0];L_x^2)\cap X^{T_0}_{b,0}\right)}\leq C(\|u_0\|_{L_x^2},R,T_0,\|k_1\|_{L_x^\infty L_y^2},a,b,c).
	$$
	
	Furthermore, according to Lemma \ref{lem:kdv_bilinear_esti} and Lemma \ref{lem:truncated_equation}, we have
	$$ 
	\begin{aligned}
	&\mathbb{E}\left\|u_{\delta,R}\right\|^2_{C([0,T];H_x^s)\cap X^{T}_{b,s} }\\
	\leq &\|u_0\|^2_{H_x^s}+(R^2+R^4)\left(\left(T\right)^{2-2(a+c)}+(T)^{1/2-a}\right)\mathbb{E}\|u_{\delta,R}\|^2_{X_{b,s}^{T}}\\
	&+(R^2T+R^4T^{1/2})\mathbb{E}\|u_{\delta,R}\|^2_{X_{b,s}^{T}} +(T+T^{2b})\|J_x^sk_1\|^{2}_{L_x^\infty L_y^2}\mathbb{E}\|u_{\delta,R}\|^{2}_{X^{T}_{b,s}},
	\end{aligned}
	$$
	for any $T>0$, similarly like the case $s=0$.
	
	Therefore, choosing $ T^*=T^*(R,\|J_x^sk_1\|_{L_x^\infty L_y^2},a,b,c)>0$ sufficiently small and dividing $[0,T_0]$ into finite numbers of intervals shorter than $T^*$,  
    we can prove 
	$$
	\mathbb{E}\|u_{\delta,R}\|^2_{C([0,T_0];H_x^s)\cap X^{T_0}_{b,s}}\lesssim C(\|u_0\|_{H_x^s},R,T_0,\|J_x^sk_1\|_{L_x^\infty L_y^2},a,b,c).
	$$
\end{proof}

\begin{rem}\label{rem:explain_local-wellposedness}
Since $ u_{\delta,R}=u_{\delta}$ in $ [0,\tau_{\delta,R,s}]$ almost surely, we have 
$$
u_{\delta} \in L^2\left(\Omega; X^{\tau_{\delta,R,s}}_{b,s}\cap C([0,\tau_{\delta,R,s}];H_x^s)\right).
$$ 
According to the proof and conditions of Lemma \ref{lem:local_well-posedness_kdv}, it is also easy to know 
$$
u_\delta\in C([0,\tau_{\delta,R,R_1,s}]; H_x^s) \cap X_{b,s}^{\tau_{\delta,R,R_1,s}},\  \forall R,R_1>0, \text{a.s.} \mathbb{P}. 
$$

What's more, 
since for any $ \varepsilon>0$, there exists a $ R_1(\varepsilon) >0$ such that $\mathbb{P}\left(\|v_\delta\|_{L_{T_0}^4H_x^s}>R_1\right)<\varepsilon$.
Thus, for any $R>0$, we can prove 
$$
u_\delta\in C([0,\tau_{\delta,R}]; H_x^s) \cap X_{b,s}^{\tau_{\delta,R}},\ \text{a.s.} \mathbb{P}.
$$
\end{rem}

To extend the solution, we need the following priori estimate.

\begin{lemma}\label{lem:priori esimate}
	Under the conditions
	$$u_0\in H_x^s,\ \|J_x^s k_1\|_{L_x^\infty L_y^2}<\infty,\ v_0\in H_x^s,\ \|J_x^sk_2\|_{L_x^\infty L_y^2}<\infty,\ \forall s\geq 0,
	$$
    we have 
	\begin{equation}\label{priori_esti_L_x^2}
		\begin{aligned}
		\mathbb{E}\|u_{\delta}\|_{ C([0,\tau_{\delta,R}];L_x^2)}^2\lesssim C(T_0,\|u_0\|_{L_x^2},\|k_1\|_{L_x^\infty L_y^2}),
		\end{aligned}
	\end{equation}
for any $R>0$, $\delta\in(0,1)$.
\end{lemma}
\begin{proof}
The proof of the case $s\in[0,2)$ will be finished through three steps. In the first step, we will introduce a series of approximation equations that the KdV type equation and dissipative  equation are both smoothed to overcome the regularity problems in using the $L_x^2$ conservation law of $u_\delta$.
In the second step, we will introduce a series of auxiliary equations, which only  their KdV type equations are smoothed. By a convergence approach, we can get a  $L_\omega^2C([0,\tau_{\delta,R,R_1}^{k'}];L_x^2)$ type prior estimate for these KdV type equations.
Finally, by a similar argument in the second step, we will illustrate \eqref{priori_esti_L_x^2}.

{\bf Step 1.}: 
Recall that  $\theta\in C_c^\infty(\mathbb{R})$, $\text{supp}\  \theta\subset[-2,2]$, $\theta|_{[-1,1]}=1$. The smooth frequency cut-off operator can be defined by 
$
\mathscr{F}^{-1}[\theta(\xi/n)\mathscr{F}[\cdot]], \ \forall n\in\mathbb{N}^+.
$
Let us introduce the following smooth initial values and noise kernels
$$
u_0^{(n)}:=\mathscr{F}^{-1}[\theta(\xi/n)\mathscr{F}[u_0]],\ 
k_1^{(n)}:=\mathscr{F}^{-1}[\theta(\xi/n)\mathscr{F}[k_1]],
$$
$$
\Phi_1^{(n)}f(x):=\int_{\mathbb{R}}k_1^{(n)}(x,y)f(y)dy.
$$

It is clear that $u_0^{(n)}\in H_x^{s_1},\ \forall s_1\geq 0$. What's more, if $k_1\in L_x^\infty L_y^2$, then by the Minkowski inequality, we have   
$$
\begin{aligned}
	\left\|J^{s_1}_x k_1^{(n)}\right\|_{L_x^\infty L_y^2}=&\left\| \mathscr{F}^{-1}[\theta(\xi/n)\left<\xi\right>^{s_1}\cdot\mathscr{F}[k_1](\xi) ]\right\|_{L_x^\infty L_y^2}\\
	\lesssim&\|\mathscr{F}^{-1}[\theta(\xi/n)\left<\xi\right>^{s_1}]\|_{L_x^1}\cdot\|k_1\|_{L_x^\infty L_y^2}\\
	\lesssim&n^{s_1}\|k_1\|_{L_x^\infty L_y^2}.
\end{aligned}
$$ 
for any $s_1\geq 0$.
The similar definitions and properties will also be gotten for $v_0^{(n)},\ k_2^{(n)},\ \Phi_2^{(n)}f(x) $.

Therefore, according to Lemma \ref{lem:well-posedness_dissipative} and Lemma \ref{lem:local_well-posedness_kdv}, for any $n\in\mathbb{N}^+$, the solution $\left(u_{\delta,R}^{(n)} ,v_{\delta,R_1}^{(n)}\right)$ of the equation 
$$ 
\begin{aligned}
	u_{\delta,R}^{(n)}=&U(t)u_0^{(n)}-\int_{0}^{t} U(t-r) \left(\frac{1}{2}\partial_{x}\left(B_R\left(u_{\delta,R}^{(n)}\right)\right)^2 +u^{(n)}_{\delta,R}\left(v^{(n)}_{\delta,R_1}\right)^2\right)dr\\
	&+\int_{0}^t U(t-r) \left(u_{\delta,R}^{(n)}\Phi^{(n)}_1 dW_r^1\right),\\
	v_{\delta}^{(n)}=&S(t)v_0^{(n)}-\delta^{-1}\int_{0}^{t} S(t-r) v_\delta^{(n)} dr+\delta^{-1/2}\int_{0}^{t} S(t-r) \left(v_\delta^{(n)} \Phi_2^{(n)} dW_r^2\right),
\end{aligned}
$$
is well-posed in 
$$ L^2\left(\Omega; X^{\tau^{n}_{\delta,R,R_1,2}}_{b,2}\cap C\left([0,\tau^n_{\delta,R,R_1,2}];H_x^2\right)\right)\times L^4\left(\Omega; L_{T_0}^4H_x^2\right) .$$

Here, we define
$$v^{(n)}_{\delta,R_1}=\theta\left(\|v^{(n)}_\delta\|_{L_t^4H_x^2}/R_1\right)v^{(n)}_\delta,\ \tau_{\delta,R}^{n}=\inf\{t\geq 0: \|u^{(n)}_{\delta,R}\|_{X^t_{b,0}}>R \}\wedge\ T_0,$$
$$ 
\tau_{\delta,R,R_1,2}^{n}=\inf\{t\geq 0: \|u^{(n)}_{\delta,R}\|_{X^t_{b,0}}>R \}\wedge\inf\{t\geq 0: \|v^{(n)}_{\delta}\|_{L_t^4H_x^2}>R_1 \}\wedge\ T_0,
$$
almost surely.

Thus, since there is no regularity problem anymore, by the It\^{o} formula, we have 
$$
\begin{aligned}
\left(u_{\delta,R}^{(n)},u_{\delta,R}^{(n)}\right)^{m}\leq&\left(u_0^{(n)},u_0^{(n)}\right)^m+2m\int_{0}^\sigma\left(u_{\delta,R}^{(n)},u_{\delta,R}^{(n)}\right)^{m-1}\left(u_{\delta,R}^{(n)}, u_{\delta,R}^{(n)}\Phi_1^{(n)}dW_r^{1}\right)\\
&+(2m^2-m)\int_{0}^\sigma \left(u_{\delta,R}^{(n)},u_{\delta,R}^{(n)}\right)^{m-1}\sum_{j=0}^\infty \left\|u_{\delta,R}^{(n)}\Phi_1^{(n)}e_j\right\|^2_{L_x^2}dr,
\end{aligned}
$$
for any $m\in\mathbb{N}^+$, where $0\leq\sigma\leq \tau_{\delta,R,R_1,2}^{n},\ \text{a.s.} \mathbb{P}$ is any arbitrary stopping time.

Hence, by the BDG inequality and the Gronwall inequality, we have the prior estimate  
\begin{equation}\label{priori_esti_L_x^2_approximation}
\begin{aligned}
	\mathbb{E}\sup_{t\in[0,\tau_{\delta,R,R_1,2}^{n}]}\left(u_{\delta,R}^{(n)},u_{\delta,R}^{(n)}\right)^{m}\leq& C(T_0,m,\|u_0\|_{L_x^2},\|k_1\|_{L_x^\infty L_y^2}),\ \forall m\in\mathbb{N}^+,
\end{aligned}
\end{equation}
which is independent to $ R,R_1$. 

Furthermore, since for any fixed $n\in\mathbb{N}^+$,  $\lim_{R_1\uparrow \infty}\inf\{t\geq 0: \|v^{(n)}_{\delta}\|_{L_t^4H_x^2}>R_1 \}\wedge T_0=T_0$ almost surely(Lemma \ref{lem:well-posedness_dissipative}) and $u^{(n)}_{\delta}(t,x)=u^{(n)}_{\delta,R}(t,x),\ t\in[0,\tau_{\delta,R,R_1,2}^{n}]$, by the monotone convergence theorem, \eqref{priori_esti_L_x^2_approximation} implies  
\begin{equation}\label{priori_esti_L_x^2_approximation_final}
	\begin{aligned}
		\mathbb{E}\sup_{t\in[0,\tau_{\delta,R}^{n}]}\left(u_{\delta}^{(n)},u_{\delta}^{(n)}\right)^{m}\leq& C(T_0,m,\|u_0\|_{L_x^2},\|k_1\|_{L_x^\infty L_y^2}),\ \forall m\in\mathbb{N}^+.
	\end{aligned}
\end{equation}

The above result \eqref{priori_esti_L_x^2_approximation_final} is also right for \eqref{moskdv} with different smoothing indexes $n_1,n_2$.
Namely, we have
\begin{equation}\label{priori_esti_L_x^2_approximation_final_different}
	\begin{aligned}
		\mathbb{E}\sup_{t\in[0,\tau_{\delta,R}^{n_1}]}\left(u_{\delta}^{(n_1,n_2)},u_{\delta}^{(n_1,n_2)}\right)^{m}\leq& C(T_0,m,\|u_0\|_{L_x^2},\|k_1\|_{L_x^\infty L_y^2}),\ \forall m\in\mathbb{N}^+,
	\end{aligned}
\end{equation}
where $u_{\delta}^{(n_1,n_2)}$ is solution of equation with initial values $u_0^{(n_1)},v_0^{(n_2)} $ and noises $\Phi_1^{(n_1)},\Phi_2^{(n_2)} $.

{\bf Step 2.:} Let us consider auxiliary equations, which only smooth the KdV type equation with any fixed index $k\in\mathbb{N}^+$
$$
\begin{aligned}
	u_{\delta,R}^{(k)'}=&U(t)u_0^{(k)}-\int_{0}^{t} U(t-r) \left(\frac{1}{2}\partial_{x}\left(B_R\left(u_{\delta,R}^{(k)'}\right)\right)^2 +u^{(k)'}_{\delta,R}\left(v_{\delta,R_1}\right)^2\right)dr\\
	&+\int_{0}^t U(t-r) \left(u_{\delta,R}^{(k)'}\Phi^{(k)}_1 dW_r^1\right),\\
	v_{\delta}=&S(t)v_0-\delta^{-1}\int_{0}^{t} S(t-r) v_\delta dr+\delta^{-1/2}\int_{0}^{t} S(t-r) \left(v_\delta \Phi_2 dW_r^2\right),
\end{aligned}
$$
We need the following stopping times and notations:
$$ 
\tau_{\delta,R,R_1}^{k'}:=\inf\{t\geq 0: \|u^{(k)'}_{\delta,R}\|_{X^t_{b,0}}>R \}\wedge\sigma_{\delta,R_1}(\omega)\wedge\ T_0,\ \mathbb{P}\text{-a.s.}$$
$$
\sigma^{n}_{\delta,R_1}(\omega)=\inf\{t\geq 0:\|v^{(n)}_{\delta}\|_{L_t^4L_x^2}\geq R_1\} \wedge T_0,\ \mathbb{P}\text{-a.s.}
$$
$$\tilde{\tau}^{n'}_{\delta,R,R_1}:=\tau_{\delta,R,R_1}^{k'}\wedge\tau_{\delta,R}^{k}\wedge \sigma^{n}_{\delta,R_1}, \ \mathbb{P}\text{-a.s.}
$$
$$
\tilde{u}^{(n)'}_{\delta,R}:=u^{(k)'}_{\delta,R}-u_{\delta,R}^{(k,n)},\ \tilde{v}^{(n)}_{\delta,R_1}:= v_{\delta,R_1}-v^{(n)}_{\delta,R_1}.
$$

Because of
\begin{equation*}
	\begin{aligned}
		v_{\delta}^{(n)}-v_\delta=&S(t)P_{>n}v_0-\delta^{-1}\int_{0}^{t} S(t-r)\left( 	v_{\delta}^{(n)}-v_\delta\right) dr\\
		&+\delta^{-1/2}\int_{0}^{t} S(t-r) \left(v_\delta \left(\Phi_2^{(n)}-\Phi_2\right) dW_r^2\right)\\
		&+\delta^{-1/2}\int_{0}^{t} S(t-r) \left(\left(v_\delta^{(n)}-v_\delta\right) \Phi_2^{(n)} dW_r^2\right),
	\end{aligned}
\end{equation*}
simple calculations show 
\begin{equation*}
	\begin{aligned}
		\left\|v_{\delta}^{(n)}-v_\delta\right\|^4_{L^4_\omega L_{T_0}^4 L_x^2}\lesssim&\left\|P_{>n}v_0\right\|^4_{L^4_\omega L_{T_0}^4 L_x^2}+\delta^{-4}T_0^3\int_{0}^{T_0} \left\|v_{\delta}^{(n)}-v_\delta\right\|^4_{L^4_\omega L_{t}^4 L_x^2} dt\\
		&+\delta^{-2}\|P_{>n}k_2\|^4_{L_x^\infty L_y^2}T^2_0\left\|v_\delta\right\|^4_{L^4_\omega L_{T_0}^4 L_x^2}\\
		&+\delta^{-2}\|k_2\|^4_{L_x^\infty L_y^2}T_0\int_{0}^{T_0} \left\|v_{\delta}^{(n)}-v_\delta\right\|^4_{L^4_\omega L_{t}^4 L_x^2}dt.
	\end{aligned}
\end{equation*}
Therefore,  by the Gronwall inequality  we have 
\begin{equation}\label{v_delta^n_convergence_result}
	\lim_{n\uparrow+\infty}\left\|v_{\delta}^{(n)}-v_\delta\right\|_{L^4_\omega L_{T_0}^4 L_x^2}=0.
\end{equation}

What's more, calculations show that for any $t\in[0,T_0]$ we have 
$$
\begin{aligned}
	&\tilde{u}^{(n)'}_{\delta,R}(t\wedge \tilde{\tau}^{n'}_{\delta,R,R_1})\\
	=&
	-\frac{1}{2}\int_{0}^{t\wedge \tilde{\tau}^{n'}_{\delta,R,R_1}}U(t\wedge \tilde{\tau}^{n'}_{\delta,R,R_1}-r)\partial_x\left[\tilde{u}_{\delta,R}^{(n)'}\left(u^{(k)'}_{\delta,R}+u^{(n)}_{\delta,R}\right)\right]dr\\
	&-\int_{0}^{t\wedge \tilde{\tau}^{n'}_{\delta,R,R_1}}U(t\wedge \tilde{\tau}^{n'}_{\delta,R,R_1}-r)\left(\tilde{u}^{(n)'}_{\delta,R}v^2_{\delta,R_1}\right)dr\\
	&-\int_{0}^{t\wedge \tilde{\tau}^{n'}_{\delta,R,R_1}}U(t\wedge \tilde{\tau}^{n'}_{\delta,R,R_1}-r)\left(u^{(k,n)}_{\delta,R}\tilde{v}^{(n)}_{\delta,R_1}\left(v^{(n)}_{\delta,R_1}+v_{\delta,R_1}\right)\right)dr\\
	&+\int_{0}^{t\wedge \tilde{\tau}^{n'}_{\delta,R,R_1}}U(t\wedge \tilde{\tau}^{n'}_{\delta,R,R_1}-r)\left(\tilde{u}^{(n)'}_{\delta,R}\Phi_1^{(k)}dW_r^{1}\right),\   \text{a.s.}\mathbb{P}.
\end{aligned}
$$

Therefore, we have 
$$
\begin{aligned}
	&\left\|\tilde{u}^{(n)'}_{\delta,R}\right\|_{C\left(\left[0,T\wedge \tilde{\tau}^{n'}_{\delta,R,R_1}\right]; L_x^2\right)\cap X^{T\wedge \tilde{\tau}^{n'}_{\delta,R,R_1}}_{b,0}}\\
	\lesssim&R\left(T^{1-(a+c)}+T^{1/4-a/2}+T^{1/2}\right)\left\|\tilde{u}^{(n)'}_{\delta,R}\right\|_{X^{T\wedge \tilde{\tau}^{n'}_{\delta,R,R_1}}_{b,0}}\\
	&+R_1^2\left(T^{1-(a+c)}+T^{1/4-a/2}+T^{1/4}\right)\left\|\tilde{u}^{(n)'}_{\delta,R}\right\|_{X^{T\wedge \tilde{\tau}^{n'}_{\delta,R,R_1}}_{b,0}}\\
	&+RR_1\left(T^{1-(a+c)}+T^{1/4-a/2}+T^{1/4}\right)\left\|\tilde{v}^{(n)}_{\delta,R_1}\right\|_{L_{T\wedge \tilde{\tau}^{n'}_{\delta,R,R_1}}^4L_x^2}\\
	&+\left\|\int_{0}^{\cdot}U(\cdot-r)\left(\tilde{u}^{(n)'}_{\delta,R}\Phi_1^{(k)}dW_r^{1}\right) \right\|_{C\left(\left[0,T\wedge \tilde{\tau}^{n'}_{\delta,R,R_1}\right];L_x^2\right)\cap X^{T\wedge \tilde{\tau}^{n'}_{\delta,R,R_1}}_{b,0}},
\end{aligned}
$$
almost surely, which implies
$$
\begin{aligned}
	&\mathbb{E}\left\|\tilde{u}^{(n)'}_{\delta,R}\right\|_{C\left(\left[0,T\wedge \tilde{\tau}^{n'}_{\delta,R,R_1}\right]; L_x^2\right)\cap X^{T\wedge \tilde{\tau}^{n'}_{\delta,R,R_1}}_{b,0}}^2\\
	\lesssim &R^2\left(T^{1-(a+c)}+T^{1/4-a/2}+T^{1/2}\right) ^2\mathbb{E}\left\|\tilde{u}^{(n)'}_{\delta,R}\right\|_{C\left(\left[0,T\wedge \tilde{\tau}^{n'}_{\delta,R,R_1}\right]; L_x^2\right)\cap X^{T\wedge \tilde{\tau}^{n'}_{\delta,R,R_1}}_{b,0}}^2\\
	&+R_1^4\left(T^{1-(a+c)}+T^{1/4-a/2}+T^{1/4}\right)^2\mathbb{E}\left\|\tilde{u}^{(n)'}_{\delta,R}\right\|^2_{X^{T\wedge \tilde{\tau}^{n'}_{\delta,R,R_1}}_{b,0}}\\
	&+R^2R^2_1\left(T^{1-(a+c)}+T^{1/4-a/2}+T^{1/4}\right)^2\mathbb{E}\left\|\tilde{v}^{(n)}_{\delta,R_1}\right\|^2_{L_{T\wedge \tilde{\tau}^{n'}_{\delta,R,R_1}}^4L_x^2}\\
	&+(T+T^{2b})\left\|k_1\right\|^2_{L_x^\infty L_y^2}\mathbb{E}\left\|\tilde{u}^{(n)'}_{\delta,R}\right\|_{X^{T\wedge\tilde{\tau}^{n'}_{\delta,R,R_1}}_{b,0}}^2.
\end{aligned}
$$

Hence, if we choose $T^*=T^*(R,R_1,\|u_0\|_{L_x^2},\|k_1\|_{L_x^\infty L_y^2},a,b,c)>0$ sufficiently small, then we have 
$$
\begin{aligned}
&\lim_{n\uparrow+\infty}\mathbb{E}\left\|\tilde{u}^{(n)'}_{\delta,R}\right\|_{C\left(\left[0,T^*\wedge \tilde{\tau}^{n'}_{\delta,R,R_1}\right]; L_x^2\right)\cap X^{T^*\wedge \tilde{\tau}^{n'}_{\delta,R,R_1}}_{b,0}}^2\\
\lesssim&	\lim_{n\uparrow+\infty}\mathbb{E}\left\|\tilde{v}^{(n)}_{\delta,R_1}\right\|^2_{L_{T^*\wedge \tilde{\tau}^{n'}_{\delta,R,R_1}}^4L_x^2}=0,
\end{aligned}
$$
by \eqref{v_delta^n_convergence_result}.
Then, by dividing  $\left[0,T_0\wedge \tilde{\tau}^{n'}_{\delta,R,R_1}\right]$ into finite numbers of intervals $\left[0,T^*\wedge \tilde{\tau}^{n'}_{\delta,R,R_1}\right]$, $[T^*\wedge \tilde{\tau}^{n'}_{\delta,R,R_1}, 2T^*\wedge \tilde{\tau}^{n'}_{\delta,R,R_1}]$,...  not longer than $ T^*$ and using arguments like \eqref{v_delta^n_convergence_result}, we can prove
\begin{equation}\label{aux_convergence_L_2}
	\lim_{n\uparrow+\infty}\mathbb{E}\left\|\tilde{u}^{(n)'}_{\delta,R}\right\|_{C\left(\left[0,T_0\wedge \tilde{\tau}^{n'}_{\delta,R,R_1}\right]; L_x^2\right)\cap X^{T_0\wedge \tilde{\tau}^{n'}_{\delta,R,R_1}}_{b,0}}^2=0.
\end{equation}

By \eqref{v_delta^n_convergence_result}, we have there exists a subsequence of $\mathbb{N}$ still denoted by $ \{n\} $ such that $	\lim_{n\uparrow+\infty}\left\|v_{\delta}^{(n)}-v_\delta\right\|_{L_{T_0}^4 L_x^2}=0 $ almost surely, which implies  \begin{equation}\label{stopping_times_2_convergence}
\lim_{n\uparrow+\infty}\sigma_{\delta,R_1}\wedge \sigma^{n}_{\delta,R_1}=\sigma_{\delta,R_1},\ \mathbb{P}\text{-a.s.}
\end{equation}

What's more, \eqref{aux_convergence_L_2} and \eqref{stopping_times_2_convergence} imply that for any $l\in\mathbb{N}^+$, there exists a $B_l\in\mathcal{F}_T$, $n(l)\in \mathbb{N}^+$ such that $\mathbb{P}(B_l)\leq 2^{-l}$ and
$$ \left\|\tilde{u}^{(n(l))'}_{\delta,R}\right\|_{ X^{\tilde{\tau}_{\delta,R,R_1,2}^{(n(l))}}_{b,0}}<R/3,\ \left\|v_{\delta}^{(n)}-v_\delta\right\|_{L_{T_0}^4 L_x^2},\  \forall  \omega\in B_l^c.
$$

Therefore, we have 
$$
\begin{aligned}
	&\mathbb{E}	\left\|u_{\delta}^{(k)'}\cdot\chi_{B_l^c}\right\|_{C\left(\left[0, \tau_{\delta,R/3,R_1/3}^{k'}\right];L_x^2\right)}^2=\mathbb{E}	\left\|u_{\delta,R}^{(k)'}\cdot\chi_{B_l^c}\right\|_{C\left(\left[0, \tau_{\delta,R/3,R_1/3}^{k'}\right];L_x^2\right)}^2\\
	\lesssim&\mathbb{E}\left\|u_{\delta, R}^{(k,n(l))}\cdot\chi_{B_l^c}\right\|^2_{C\left(\left[0, \tilde{\tau}_{\delta,R,R_1}^{n(l)'}\right]; L_x^2\right)}+\mathbb{E}\left\|\tilde{u}_{\delta, R}^{(n(l)')}\cdot\chi_{B_l^c}\right\|^2_{C\left(\left[0, \tilde{\tau}_{\delta,R,R_1}^{n(l)'}\right]; L_x^2\right)}.
\end{aligned}
$$

What's more, since $\{\chi_{B_l^c}\}_{l\in\mathbb{N}^+}$ converge to $1$ in probability, there exists a subsequence still denoted by $ \{B_l^c\}_{l\in\mathbb{N}^+}$ converges to 
$1$ almost surely.

Hence, by the Fatou Lemma and \eqref{priori_esti_L_x^2_approximation_final_different}, we have 
\begin{equation*}
	\begin{aligned}
		&\mathbb{E}	\left\|u_{\delta}^{(k)'}\right\|_{C\left(\left[0, \tau_{\delta,R/3,R_1/3}^{k'}\right];L_x^2\right)}^2\\
		\lesssim&\liminf_{l\uparrow\infty}\mathbb{E}	\left\|u_{\delta}^{(k)'}\cdot\chi_{B_l^c}\right\|_{C\left(\left[0, \tau_{\delta,R/3,R_1/3}^{k'}\right];L_x^2\right)}^2\\
		\lesssim&\lim_{l\uparrow\infty}\mathbb{E}\left\|u_{\delta, R}^{(k,n(l))}\right\|^2_{C\left(\left[0, \tilde{\tau}_{\delta,R,R_1}^{n(l)'}\right]; L_x^2\right)}+\lim_{l\uparrow\infty}\mathbb{E}\left\|\tilde{u}_{\delta, R}^{(n(l)')}\right\|^2_{C\left(\left[0, \tilde{\tau}_{\delta,R,R_1}^{n(l)'}\right]; L_x^2\right)}\\
		\lesssim&C(T_0,\|u_0\|_{L_x^2},\|k_1\|_{L_x^\infty L_y^2}),
	\end{aligned}
\end{equation*}
which implies 
\begin{equation}\label{auxiliary_stopping_time_result_dot}
	\mathbb{E}	\left\|u_{\delta}^{(k)'}\right\|_{C\left(\left[0, \tau_{\delta,R}^{k'}\right];L_x^2\right)}^2	\lesssim C(T_0,\|u_0\|_{L_x^2},\|k_1\|_{L_x^\infty L_y^2}),
\end{equation}
because of $\lim_{R_1\uparrow \infty}\sigma_{\delta,R_1}(\omega)\wedge\ T_0=T_0,\ \mathbb{P}\text{-a.s.}$

{\bf Step 3.:} Let us set 
$$ \tilde{\tau}^{k}_{\delta,R,R_1}:=\tau_{\delta,R,R_1}\wedge\tau_{\delta,R,R_1}^{k'},\
\tilde{u}^{(k)}_{\delta,R}:=u_{\delta,R}-u_{\delta,R}^{(k')},\ \mathbb{P}\text{-a.s.} 
$$

According to the fact  

Moreover, simple calculations show that for any $t\in[0,T_0]$ 
$$
\begin{aligned}
	&\tilde{u}^{(k)}_{\delta,R}(t\wedge \tilde{\tau}^{k}_{\delta,R,R_1})\\
	=&U(t\wedge \tilde{\tau}^{k}_{\delta,R,R_1})\left(u_0-u_0^{(k)}\right)\\
	&
	-\frac{1}{2}\int_{0}^{t\wedge \tilde{\tau}^{k}_{\delta,R,R_1}}U(t\wedge \tilde{\tau}^{k}_{\delta,R,R_1}-r)\partial_x\left[\tilde{u}_{\delta,R}^{(k)}(u_{\delta,R}+u^{(k)'}_{\delta,R})\right]dr\\
	&-\int_{0}^{t\wedge \tilde{\tau}^{k}_{\delta,R,R_1}}U(t\wedge \tilde{\tau}^{k}_{\delta,R,R_1}-r)\left(\tilde{u}^{(k)}_{\delta,R}v^2_{\delta,R_1}\right)dr\\
	&+\int_{0}^{t\wedge \tilde{\tau}^{k}_{\delta,R,R_1}}U(t\wedge \tilde{\tau}^{k}_{\delta,R,R_1}-r)\left(\tilde{u}^{(k)}_{\delta,R}\Phi_1dW_r^{1}\right)\\
	&+\int_{0}^{t\wedge \tilde{\tau}^{k}_{\delta,R,R_1}}U(t\wedge \tilde{\tau}^{k}_{\delta,R,R_1}-r)\left(u_{\delta,R}\left(\Phi_1-\Phi_1^{(k)}\right)dW_r^{1}\right),\   \text{a.s.}\mathbb{P}.
\end{aligned}
$$

Thus, we have 
$$
\begin{aligned}
	&\mathbb{E}	\left\|\tilde{u}^{(k)}_{\delta,R}\right\|_{C\left(\left[0,T\wedge \tilde{\tau}^{k}_{\delta,R,R_1}\right];L_x^2\right)\cap X^{T\wedge\tilde{\tau}^{k}_{\delta,R,R_1}}_{b,0}}^2\\
	\lesssim&C(T)\left\|u_0-u_0^{(k)}\right\|_{L_x^2}^2+R^2\left(T^{1-(a+c)}+T^{1/4-a/2}+T^{1/2}\right)^2\mathbb{E}\left\|\tilde{u}^{(k)}_{\delta,R}\right\|_{X^{T\wedge\tilde{\tau}^{k}_{\delta,R,R_1}}_{b,0}}^2\\
	&+R^4_1\left(T^{1-(a+c)}+T^{1/4-a/2}+T^{1/4}\right)^2\mathbb{E}\left\|\tilde{u}^{(k)}_{\delta,R}\right\|_{X^{T\wedge\tilde{\tau}^{k}_{\delta,R,R_1}}_{b,0}}^2\\
	&+(T+T^{2b})\left\|k_1\right\|^2_{L_x^\infty L_y^2}\mathbb{E}\left\|\tilde{u}^{(k)}_{\delta,R}\right\|_{X^{T\wedge\tilde{\tau}^{k}_{\delta,R,R_1}}_{b,0}}^2\\
	&+(T^{1-2b}+1)\mathbb{E}\left\|u_{\delta,R}\left(\Phi_1-\Phi_1^{(k)}\right)e_j\right\|_{L^2_{ T\wedge\tilde{\tau}^{k}_{\delta,R,R_1}}l^2_j L_x^2}^2.
\end{aligned}
$$

Hence, by the argument in Step 2 and Lemma 5.3 in \cite{kdvas}, we have 
$$
\lim_{k\uparrow+\infty}\mathbb{E}	\left\|\tilde{u}^{(k)}_{\delta,R}\right\|_{C\left(\left[0,T_0\wedge \tilde{\tau}^{k}_{\delta,R,R_1}\right];L_x^2\right)\cap X^{T_0\wedge\tilde{\tau}^{k}_{\delta,R,R_1}}_{b,0}}^2=0,
$$
which implies 
\begin{equation*}
	\begin{aligned}
		\mathbb{E}\|u_{\delta}\|_{ C([0,\tau_{\delta,R}];L_x^2)}^2\lesssim C(T_0,\|u_0\|_{L_x^2},\|k_1\|_{L_x^\infty L_y^2}),
	\end{aligned}
\end{equation*}
like \eqref{auxiliary_stopping_time_result_dot}.
\end{proof}

Recall Remark \ref{rem:explain_local-wellposedness}, we have $u_\delta\in C([0,\tau_{\delta,R}]; H_x^s) \cap X_{b,s}^{\tau_{\delta,R}} $ almost surely. Now, we extend the solution to  $[0,T_0]$ almost surely. Let us define $\tau^*_\delta:=\lim_{R\uparrow \infty}\tau_{\delta,R},\ \text{a.s.} \mathbb{P}$

\begin{prop}\label{prop:global_well-posedness}
For any $T_0>0$, under the conditions of Lemma \ref{lem:priori esimate}, we have $u_\delta\in C([0,T_0];H_x^s)\cap X^{T_0}_{b,s},\ \forall b\in(5/12,1/2),\ \text{a.s.} \mathbb{P}$.
\end{prop}
\begin{proof}
{\bf The case $s=0$.:}	
We argue this case by a contradiction.

It is sufficient to prove $ \tau_{\delta}^*=T_0$ almost surely, for any $T_0>0$. If not, there exists  a $T_0>0$ and  $\lambda >0,\ E:=\{\omega:\tau^*_{\delta}<T_0\}\in \mathcal{F}_{T_0}$ such that $\mathbb{P}(E)=\lambda>0$.
	
	According to $\lim_{R\uparrow \infty}\tau_{\delta,R} =\tau^*_\delta\ \text{a.s.} \mathbb{P}$, for any $\varepsilon>0$, there exists a $R = R(\varepsilon)>0$ such that 
	\begin{equation}\label{contradiction_1}
	\mathbb{P}(\{\tau^*_{\delta}-\tau_{\delta,R}>\varepsilon\})<\lambda/4.
	\end{equation}

	Let us consider the equation starting from $\tau_{\delta,R}$ adapted to $ \mathcal{F}_{t+\tau_{\delta,R}}$
\begin{equation}\label{equation_star_from_stopping_times}
	\begin{aligned}
	\mu_{\delta,\tilde{R}}(t)=&U(t)u_{\delta,R}(\tau_{\delta,R})-\int_{0}^{t} U(t-r)(B_{\tilde{R}}(\mu_{\delta,\tilde{R}}))_xdr\\
	&-\int_{0}^{t} U(t-r)(\mu_{\delta,\tilde{R}}\nu_{\delta,\tilde{R}}^2)dr+\int_{0}^{t}U(t-r) (\mu_{\delta,\tilde{R}}\Phi_1 dW_{r+\tau_{\delta,R}}),
\end{aligned}
\end{equation}
	where $\tilde{R}>0$ is undetermined and $\nu_{\delta,\tilde{R}}(t)=\theta_{\tilde{R}}(\|v_{\delta}(\tau_{\delta,R}+\cdot)\|_{L_t^4L_x^2})v_{\delta}(\tau_{\delta,R}+\cdot)$. What's more, according to the strong Markov property of the Brownian motion,  $W_{t'+\tau_{\delta,R}}-W_{\tau_{\delta,R}}$ is also a Brownian motion.
	
Let us set    
$
\iota_{\delta,\tilde{R}}:=\inf\{t\geq0:\|\nu_{\delta}\|_{L_t^4L_x^2}\geq\tilde{R}\},\ \text{a.s.}\mathbb{P}
$.
According to Lemma \ref{lem:well-posedness_dissipative}, we can choose $\tilde{R}>C(\|v_0\|_{L_x^2},T_0,\|k_2\|_{L_x^\infty L_y^2})$ sufficiently large such that
\begin{equation}\label{choose_c_R_vdelta}
\mathbb{P}\left(\iota_{\delta,\tilde{R}}<T_0\right)<\lambda/4.
\end{equation}

From the proof of Lemma \ref{lem:local_well-posedness_kdv}, for any  $ \tilde{R}>0$ satisfying \eqref{choose_c_R_vdelta} and 
$$
\frac{4\mathbb{E}\|u_\delta\|_{C([0,\tau_\delta^*];L_x^2)}^2}{\tilde{R}^2}<
\lambda/4.
$$
there exists a
$\tilde{T} = \tilde{T}(\tilde{R},\|k_1\|_{L_x^\infty L_y^2})\in(0,T_0)$ sufficiently small such that 
$$
\begin{aligned}
	\|\mu_{\delta,\tilde{R}}\|_{L^2_\omega X^{\tilde{T}}_{b,0}}\leq& 2\|u_{\delta,R}(\tau_{\delta,R})\|_{L^2_\omega L^2_x}\\
	=& 2\|u_\delta\|_{L^2_\omega \left(C([0,\tau_\delta^*];L_x^2)\right)}\\
	\leq& C\left(T_0,\|u_0\|_{L_x^2},\|k_1\|_{L_x^\infty L_y^2}\right).
\end{aligned}
$$
	
Thus, we have  
\begin{equation}\label{extend_u_from_stopping_times}
	\begin{aligned}
	\mathbb{P}\left(\{\|\mu_{\delta,\tilde{R}}\|_{X^{\tilde{T}}_{b,0}}\geq\tilde{R}\}\right)\leq \frac{C\left(T_0,\|u_0\|_{L_x^2},\|k_1\|_{L_x^\infty L_y^2}\right)^2}{\tilde{R}^2}<
	\lambda/4.
\end{aligned}
\end{equation}

	Therefore, according to \eqref{choose_c_R_vdelta} and \eqref{extend_u_from_stopping_times}, we have 
	\begin{equation}\label{contradiction_2}
	\begin{aligned}
		&\mathbb{P}\left(\{ \tau_{\delta,R}+\tilde{T}< \tau_{\delta}^*<T_0\}\right)\\
		\geq &\mathbb{P}\left(\{\|\mu_{\delta,\tilde{R}}\|_{X^{\tilde{T}}_{b,0}}<\tilde{R}\}\cap\{v_\delta=v_{\delta,\tilde{R}},\ t\in[0,T_0]\}\cap E\right)\\
		=&\mathbb{P}(E)-\mathbb{P}\left(\{\|\mu_\delta\|_{X^{\tilde{T}}_{b,0}}\geq\tilde{R}\}\right)-\mathbb{P}\left(\{\iota_{\delta,\tilde{R}}<\tilde{T}\}\right)\\
		\geq& \lambda/2.
	\end{aligned}
	\end{equation}

Thus, if we let $\varepsilon = \tilde{T}/2$, then \eqref{contradiction_2}  will contradict to \eqref{contradiction_1}.

Hence, we have finished the proof of the case $s=0$.

    {\bf The case $s>0$.:} 
    Based on Lemma \ref{lem:embedding_result}, we have $\|J^sk_i(x,y)\|_{L_x^{\infty}L_y^2}<\infty$ implies $\|k_i(x,y)\|_{L_x^{\infty}L_y^2}<\infty $, for any $s\geq0$. Thus, the proof of  case $s=0$ implies $ \lim_{R\uparrow \infty} \tau_{\delta,R}\wedge T_0=T_0$ almost surely, for any $T_0>0$. Thus, according to Remark \ref{rem:explain_local-wellposedness} we can prove $u_\delta\in C([0,T_0];H_x^s)\cap X^{T_0}_{b,s}$, for any $T_0>0$.

\end{proof}

\begin{lemma}\label{lem:embedding_result}
	For any $ s>0$,  $\|J_x^sf(x)\|_{L_x^\infty}<\infty$ implies $\|f(x)\|_{L_x^\infty}<\infty$.
\end{lemma}
\begin{proof}
	Because the Minkowski inequality implies
	$$
	\|f\|_{L_x^\infty}=\|J_x^{-s}J_x^sf\|_{L_x^\infty}\lesssim\|\mathscr{F}^{-1}[(1+|\xi|^2)^{-s/2}](x)\|_{L_x^1}\|J_x^sf\|_{L_x^\infty},
	$$
 we only need to prove 
	$$ \|\mathscr{F}^{-1}[(1+|\xi|^2)^{-s/2}]\|_{L_x^1}<\infty.
	$$
	
	Since 
	$$
	(1+|\xi|^2)^{-s/2}=\frac{1}{\Gamma(s/2)}\int_{0}^\infty e^{-(1+|\xi|^2)u}u^{\frac{s}{2}-1} du,
	$$
	we have 
	$$
	\begin{aligned}
	\mathscr{F}^{-1}[(1+|\xi|^2)^{-s/2}](x)=&\frac{1}{\Gamma(s/2)}\int_{0}^\infty e^{-u}u^{\frac{s}{2}-1}\int_{R_\xi}  e^{ix\xi} e^{-\xi^2u}d\xi du\\
	=&\frac{1}{\Gamma(s/2)}\int_{0}^\infty e^{-u}u^{\frac{s}{2}-1} e^{-cx^2/u}du>0.
	\end{aligned}
	$$
	Thus, we have 
	$$
	\begin{aligned}
	\|\mathscr{F}^{-1}[(1+|\xi|^2)^{-s/2}]\|_{L_x^1}\lesssim&\int_{-\infty}^\infty\int_{0}^\infty e^{-u}u^{\frac{s}{2}-1} e^{-cx^2/u}dudx<\infty,
	\end{aligned}
	$$
	which finishes the proof.
\end{proof}

Now, we can complete the proof of the well-posedness of \eqref{moskdv}. 
\begin{proof}[{\bf Proof of Theorem \ref{thm:global_well-posedness}:}]
Combining Lemma \ref{lem:well-posedness_dissipative} with Proposition \ref{prop:global_well-posedness}, we have finished the proof.
\end{proof}

\section{A weak type averaging principle of stochastic KdV equations}\label{sec:averaging_principle}

In this section, we will prove the averaging principle in the sense of probability. 
According to Corollary \ref{coro:longtime_dissipative_equation},  the average equation of the initial system \eqref{moskdv} should be \eqref{average_eq}, which is a stochastic KdV equation.  The well-posedness of \eqref{average_eq} has been studied in \cite{kdvas}.

\begin{lemma}(\cite{kdvas})\label{lem:average_eq}
For any fixed  $T_0>0,\ s\geq 0$ and $\forall b\in(3/8,1/2)$, under the conditions $u_0\in H_x^s$ and $\|J_x^s k_1\|_{L_x^\infty L_y^2}<\infty$, \eqref{average_eq} is global well-posed in $ C([0,T_0];H_x^s)\cap X^{T_0}_{b,s}$ almost surely. 

\end{lemma}

The uniform tail estimate $ \mathbb{P}\left(\|u_\delta\|_{X^{T_0}_{b,s}}\geq R \right)$ with respect to $\delta$ in the next lemma will be used in the proof of Theorem \ref{thm:average_principle}. Since we only have moment estimates of $u_{\delta,R}$, thus we achieve this goal by using the property of $u$ and a convergence method.

\begin{lemma}\label{lem:probability_estimate}
Under the conditions of Theorem \ref{thm:average_principle}, for any  $T_0>0$, $\zeta\in(0,1)$, there exists a $R_0=R_0(\zeta,T_0)>0$  such that for any fixed $ R>R_0$, it holds
\begin{equation}\label{meas_esti}
	\begin{aligned}
		&\mathbb{P}\left(\|u\|_{X^{T_0}_{b,s}}\geq R \right)\leq \zeta,\\ &\mathbb{P}\left(\|u_\delta\|_{X^{T_0}_{b,s}}\geq R \right)\leq \zeta,\ \forall \delta\in(0,\delta_0),
	\end{aligned}
\end{equation}
where $\delta_0$ is a constant depending on $ R,\zeta,T_0$. 

\end{lemma}
\begin{proof}
	According to Lemma \ref{lem:average_eq}, there exists a $R^*>0$ such  that 
\begin{equation}\label{tail_esti_1}
	\mathbb{P}\left(\|u\|_{X^{T_0}_{b,s}}>2R^*/3\right)<\zeta/3.
\end{equation}
	
	For any $s\geq 0$, $\delta\in(0,1)$, $R>0$, let us set 
	$$
	\tau^s_{R}:=\inf\{t\geq0:\|u\|_{X^t_{b,s}}\geq R\}\wedge T_0,
	$$
	$$
	\tau^s_{\delta,R}:=\inf\{t\geq0:\|u_\delta\|_{X^t_{b,s}}\geq R\}\wedge T_0,
	$$
	$$
	\tilde{\tau}^s_{\delta,R}:=\tau^s_{\delta,R}\wedge\tau_{R}^s,
	$$
	almost surely.
	
	By calculations, we have 
	$$
	\begin{aligned}
		&\mathbb{E}\|u-u_{\delta}\|^2_{C([0,\tilde{\tau}^s_{\delta,R}\wedge T];H_x^s)\cap X^{\tilde{\tau}^s_{\delta,R}\wedge T}_{b,s}}\\
		\lesssim&\mathbb{E}\left\|\int_{0}^{t} U(t-r)\partial_x((u-u_\delta)(u+u_\delta))dr\right\|_{C([0,\tilde{\tau}^s_{\delta,R}\wedge T];H_x^s)\cap X^{\tilde{\tau}^s_{\delta,R}\wedge T}_{b,s}}^2\\
		&+\mathbb{E}\left\|\int_{0}^{t} U(t-r)u_\delta v_\delta^2dr\right\|_{C([0,\tilde{\tau}^s_{\delta,R}\wedge T];H_x^s)\cap X^{\tilde{\tau}^s_{\delta,R}\wedge T}_{b,s}}^2\\
		&+\mathbb{E}\left\|\int_{0}^{t} U(t-r)((u-u_{\delta})\Phi_1dW_r^{1})\right\|_{C([0,\tilde{\tau}^s_{\delta,R}\wedge T];H_x^s)\cap X^{\tilde{\tau}^s_{\delta,R}\wedge T}_{b,s}}^2,
	\end{aligned}
	$$
	for any $\delta\in(0,1)$, $T\in(0,T_0)$, $b\in(5/12,1/2)$.
	
	Thus, just like Lemma \ref{lem:local_well-posedness_kdv}, we have 
	$$
	\begin{aligned}
		&\mathbb{E}\|u-u_{\delta}\|^2_{C([0,\tilde{\tau}^s_{\delta,R}\wedge T];H_x^s)\cap X^{\tilde{\tau}^s_{\delta,R}\wedge T}_{b,s}}\\
		\lesssim& \left(T^{2-2(a+c)}+T^{1/2-a}+T\right)R^2\mathbb{E}\|u-u_{\delta}\|^2_{X^{\tilde{\tau}^s_{\delta,R}\wedge T}_{b,s}}\\
		&+\left(T^{2-2(a+c)}+T^{1/2-a}+T^{1/2}\right)R^2\mathbb{E}\|v_{\delta}\|^4_{L_T^4H_x^s}\\
		&+(T+T^{2b})\|J_x^sk_1\|^{2}_{L_x^\infty L_y^2}\mathbb{E}\|u-u_{\delta}\|^2_{X^{\tilde{\tau}^s_{\delta,R}\wedge T}_{b,s}},
	\end{aligned}
	$$
	for any $a,b\in(5/12,1/2),\ c>1/2,\ a+c<1$. Without loss of generality, we further assume $ a+2c<3/2$.

	Therefore, if we choose $R,T>0$ such that
	$$ R>R^*,\ T\lesssim \|J_x^sk_1\|^{-b^{-1}}_{L_x^\infty L_y^2},\  T\lesssim    R^{-8/(1-2a)},$$
	then we have 
	\begin{equation*}
		\begin{aligned}
			\mathbb{E}\|u-u_{\delta}\|^2_{C([0,\tilde{\tau}^s_{\delta,R}\wedge T];H_x^s)\cap X^{\tilde{\tau}^s_{\delta,R}\wedge T}_{b,s}}
			\lesssim 	R^2\mathbb{E}\|v_{\delta}\|^4_{L_{T}^4H_x^s}.
		\end{aligned}
	\end{equation*}
	
	Therefore, if we dividing $[0,T_0] $ into finite numbers of intervals $[0,T]$, $[T,2T]$,..., $[T_0-T,T_0]$, we have 
	\begin{equation}\label{preparation_for_convergence}
		\begin{aligned}
			\mathbb{E}\|u-u_{\delta}\|^2_{C([0,\tilde{\tau}^s_{\delta,R}\wedge T_0];H_x^s)\cap X^{\tilde{\tau}^s_{\delta,R}\wedge T_0}_{b,s}}
			\lesssim 	C(R,\|J_x^sk_1\|_{L_x^\infty L_y^2})\mathbb{E}\|v_{\delta}\|^4_{L_{T_0}^4H_x^s}.
		\end{aligned}
	\end{equation}	
	
	According to Lemma \ref{lem:well-posedness_dissipative}, because of $\lim_{\delta\downarrow 0}\mathbb{E}\|v_{\delta}\|^4_{L_{T_0}^4H_x^s}\sim\lim_{\delta\downarrow 0}\delta= 0$,
	we can choose $\delta_{0}(\zeta,R)>0$ such that  for any $\delta<\delta_{0}$, there exist a $B_{\delta}\subset \mathcal{F}_{T_0}$ satisfying
\begin{equation}\label{tail_esti_2}
	\mathbb{P}(B_\delta)\lesssim \delta <\zeta/3,\ \Omega\backslash B_\delta=\left\{\omega: \|u-u_\delta\|_{X^{\tilde{\tau}^s_{\delta,R}\wedge T_0}_{b,s}}\leq R/3\right\}. 
\end{equation}
	
	Thus, it follows
\begin{equation}\label{tail_esti_3}
	\begin{aligned}
	\mathbb{E}\| u_{\delta}\cdot \chi_{\Omega\backslash B_\delta}(\omega)\|^2_{X^{\tau^s_{2R/3}\wedge T_0}_{b,s}}\leq&\mathbb{E}\| u_{\delta}\cdot \chi_{\Omega\backslash B_\delta}(\omega)\|^2_{X^{\tilde{\tau}^s_{\delta,R}\wedge T_0}_{b,s}} \\
	\leq&\mathbb{E}\| u-u_{\delta}\|^2_{X^{\tilde{\tau}^s_{\delta,R}\wedge T_0}_{b,s}}+\mathbb{E}\| u\|^2_{X^{\tau^s_{R}\wedge T_0}_{b,s}} \\
	\leq &2R^2,
	\end{aligned}
\end{equation}
	for any $ \delta\in(0,\delta_0).$
	
	Hence, by \eqref{tail_esti_1}, \eqref{tail_esti_2}, \eqref{tail_esti_3} and the Chebyshev inequality, we have 
\begin{equation}\label{tail_esti_4}
	\begin{aligned}
		&\mathbb{P}\left(\| u_{\delta}\|_{X_{b,s}^{T_0}}\geq (\zeta/6)^{-1/2}R\right)\\
		\leq&\mathbb{P}\left(\left\{\| u_{\delta}\|_{X_{b,s}^{T_0}}\geq (\zeta/6)^{-1/2}R\right\}\cap  \left\{\omega: \tau^s_{2R/3}(\omega ) \geq T_0\right\} \cap \{\Omega\backslash B_\delta\} \right)\\
		&+\mathbb{P}( \tau^s_{2R/3}(\omega )< T_0)+\mathbb{P}(B_\delta)\\
		\leq& \mathbb{P}\left( \| u_{\delta}\cdot \chi_{\Omega\backslash B_\delta}(\omega)\|_{X^{\tau^s_{2R/3}\wedge T_0}_{b,s}}\geq (\zeta/6)^{-1/2}R \right)+2\zeta/3\\
		\leq&\zeta, \ \forall R>R^*,\ \delta \in(0,\delta_{0}).
	\end{aligned}
	\end{equation}
	
	Combining \eqref{tail_esti_1} with \eqref{tail_esti_4}, by choosing $R_0=(\zeta/6)^{-1/2}R^*\vee R^*$, we can finish the proof.

\end{proof}

\begin{proof}[{\bf Proof of Theorem \ref{thm:average_principle}:}]
	
Firstly, according to the conditions of Theorem \ref{thm:average_principle} and Corollary \ref{coro:longtime_dissipative_equation}, \eqref{average_eq} is the so-called averaging equation of \eqref{moskdv}.

By Lemma \ref{lem:probability_estimate}, for any  $T_0,\zeta/2>0$, we can choose $R, \delta_0>0$ such that
\begin{equation}\label{final_esti}
\mathbb{P}(\{\tilde{\tau}^s_{\delta,R}< T_0 \} )<\zeta/2,\ \forall \delta\in(0,\delta_0).
\end{equation}
	
Under the restriction of Lemma \ref{lem:probability_estimate}, according to \eqref{preparation_for_convergence}, we have
the following estimate 
$$
\begin{aligned}
\mathbb{E}\|u-u_{\delta}\|^2_{C([0,\tilde{\tau}^s_{\delta,R}\wedge T_0];H_x^s)}
\lesssim C(R)\mathbb{E}\|v_{\delta}\|^4_{L_{T_0}^4H_x^s}.
\end{aligned}
$$
for any $\delta\in(0,1)$, $T\in(0,T_0)$, $b\in(5/12,1/2)$.

Therefore, by \eqref{final_esti} and  the Chebyshev inequality, there exists a $\delta_0=\delta_0(R,\zeta,\alpha)>0$
$$
\begin{aligned}
	&\mathbb{P}(\|u-u_{\delta}\|_{C([0,T_0];H_x^s)\cap X_{b,s}^{T_0}}>\alpha)\\
	\leq& \mathbb{P}(\|u-u_{\delta}\|_{C([0,T_0];H_x^s)\cap X_{b,s}^{T_0}}>\alpha\cap\{\tilde{\tau}^s_{\delta,R}\geq T_0 \} )+\mathbb{P}(\{\tilde{\tau}^s_{\delta,R}< T_0 \} )\\
	\leq&\mathbb{P}\left(\|u-u_{\delta}\|_{C([0,\tilde{\tau}^s_{\delta,R}\wedge T_0];H_x^s)\cap X_{b,s}^{\tilde{\tau}^s_{\delta,R}\wedge T_0}}>\alpha\right)+\mathbb{P}(\{\tilde{\tau}^s_{\delta,R}< T_0 \} )\\
	<&\zeta,
\end{aligned}
$$
for any $\delta\in(0,\delta_0)$.

Hence, we have finished the proof.
\end{proof}

\begin{rem}\label{rem:convergence_rate_probability}
	Different from Theorem \ref{thm:result_with_classical_workspace}, according to the proof of Lemma \ref{lem:probability_estimate}, because we don't know the relationship between the $\zeta$ and $R^*$,  we can not get a explicit relationship between  $\delta_{0}$ and $ \zeta$. Thus, we can not get a probability convergence rate with respect to $\delta$.
\end{rem}

\section{A strong type averaging principle of stochastic KdV equations}\label{sec:averaging_principle_with_moment}

In this section, we prove the averaging principle of \eqref{moskdv} in the $L_\omega^2$ sense.  We will use the well-posedness result of multiplicative stochastic KdV equation in \cite{kdv_sharper}, since \cite{kdv_sharper} is a recent refinement of \cite{kdvmuti}. Roughly speaking, the main difference between \cite{kdvmuti} and \cite{kdv_sharper} is that in \cite{kdv_sharper} we use a sharper bilinear estimate and a workspace slightly different from the deterministic case to lower the requirement  of the noise kernel form $ L_x^2\cap L_x^1$ to $L_x^2\cap L_x^{2/(1+4q)}, \ \forall 0<q\ll1$.  We also want to point out the approach to prove the global well-posedness  of \cite{kdv_sharper} is similar to that of \cite{kdvmuti}.

In this section, just like \cite{kdvmuti} and \cite{kdv_sharper}, we consider the convolution type noise kernel. Namely, we let 
$$\Phi_1f(x)=\int_{\mathbb{R}} f(y)k_1(x-y)dy.$$
This is because although we can treat general noise kernels like \eqref{auxiliary_equality},  convolution type noise enables us to clearly observe the non-regularity demands that the homogeneous low-frequency restriction of workspace places on the noise kernel.

\begin{rem}\label{rem:why_L_x2}
	According to \cite{kdvmuti} or \cite{kdv_sharper}, we know that if we want to prove the well-posedness of KdV type equation in $L^{2m}_\omega \left(C([0,T_0]; H_x^s)\cap Z_{b,s,q}^{T_0} \right)$, the $H_x^s$-conservation law of KdV type equation is needed. Thus, owing to the corresponding deterministic equation of \eqref{moskdv} only has $L_x^2$-conservation law, we can only consider the well-posedness and the averaging principle for the case $s=0$ in this part. 
\end{rem}

To prove the local well-posedness of $u_\delta$ in $L^{2m}_\omega \left(C([0,T_0]; L_x^2)\cap Z_{b,0,q}^{T_0} \right)$, the following localized  modified equation is needed 
\begin{equation}\label{cutted_off_2_equation}
	\begin{aligned}
	w_{\delta,R}(t,x)=&-\int_{0}^{t} U(t-r)\partial_{x}((\tilde{w}_{\delta,R}+U(r)u_0)^2)dr \\
	&-\int_{0}^{t} U(t-r) \left((\tilde{w}_{\delta,R}+U(r)u_0)v_{\delta,R}^2\right)dr\\
	&+\int_{0}^{t} U(t-r)\left((\tilde{w}_{\delta,R}+U(r)u_0)\Phi_1 dW^{(1)}_r\right),
	\end{aligned}
\end{equation}
for any $R>0$.
Here, we have 
$$
w_{\delta,R}=u_{\delta,R}-U(t)u_0,\   \tilde{w}_{\delta,R}(t,x)=\theta_R\left(\|w_{\delta,R}\|_{Z^t_{b,0,q}}\right)w_{\delta,R},\ w_\delta=u_\delta-U(t)u_0.
$$

Let us prove the local well-posedness of \eqref{moskdv} in  $u_{\delta,R}\in L_{\omega}^{2m}\left(Y^T_{b,0}\cap C([0,T_0];L_x^2)\right),$

\begin{lemma}\label{lem:local_well-posedness_moment_estimate}
For any fixed $T_0>0$, and any $R\in (0,+\infty)$, $ \delta\in(0,1)$,  $  q\in (0,2/19)$, $ a\in(5/12-q/3,1/2),\ b\in\left(1/2-q^2/3,1/2\right)$, $c>1/2$, $b+3a>3/2$, $a+c<1$, we have
$$
w_{\delta,R}\in L_{\omega}^{2m}\left(Z^T_{b,0,q}\cap C([0,T_0];L_x^2\cap \dot{H}_{x}^{-q})\right),\ u_{\delta,R}\in L_{\omega}^{2m}\left(Y^T_{b,0}\cap C([0,T_0];L_x^2)\right),
$$
for any $m\in\mathbb{N}^+$.
\end{lemma}
\begin{proof}

	According to Corollary \ref{coro:restricted_estimate}, Proposition 2.3 in \cite{kdv_sharper}, Theorem 1.1 in \cite{kdv_sharper}, Lemma \ref{lem:trilinear_estimate_q} and Lemma \ref{lem:stochastic_convolution_estimate_2}, for any $T\in(0,1)$, we have
	\begin{equation}\label{2_estimate_local}
	\begin{aligned}
    &\mathbb{E}\|w_{\delta,R}-w_{\delta,R}'\|^{2m}_{Z^T_{b,0,q}\cap C([0,T];L_x^2\cap \dot{H}_{x}^{-q})}\\
    \lesssim & (R+\|u_0\|_{L_x^2})^{2m}\left(T^{1-(a+c)}+T^{1/4-a/2}\right)^{2m} \mathbb{E}\|w_{\delta,R}-w_{\delta,R}'\|^{2m}_{Z^T_{b,0,q}}\\
    &+R^{2m}\left(T^{1-(a+c)}+T^{1/4-a/2}\right)^{2m}\mathbb{E}\|w_{\delta,R}-w_{\delta,R}'\|^{2m}_{Y^T_{b,0}}\\
    &+T^{2bm}\left(\|k_1\|^{2m}_{L_x^{2/(1+4q)}}+\|k_1\|_{L_x^2}^{2m}\right)\mathbb{E}\|w_{\delta,R}-w_{\delta,R}'\|^{2m}_{Y^T_{b,0}},
	\end{aligned}
	\end{equation}
under the restrictions of $a,b,c ,q $ in this lemma.  

Thus, choosing $T=T(R,\|u_0\|_{L_x^2},\|k_1\|_{L_x^2},\|k_1\|_{L_x^{2/(1+4q)}})$ sufficiently small, by a fixed point argument, we can get the well-posedness of $w_{\delta,R}$ in $$L_{\omega}^{2m}\left(Z^T_{b,0,q}\cap C([0,T];L_x^2\cap \dot{H}_{x}^{-q})\right).$$

Note that for any $\Gamma \in(0,T_0),\ t\in[0,T_0-\Gamma]$,  we have 
$$
\begin{aligned}
w_{\delta,R}(\Gamma+ t)=& U(t)u_{\delta,R}(\Gamma)-U(t+\Gamma)u_0\\
&-\int_{0}^{t} U(t-r)\partial_x(\tilde{w}_{\delta,R}(\Gamma+r)+U(\Gamma+r)u_0)^2 dr\\
&-\int_{0}^{t} U(t-r)\left((\tilde{w}_{\delta,R}(\Gamma+r)+U(\Gamma+r)u_0)v_{\delta,R}^2(\Gamma+r)\right) dr\\
&+\int_{0}^{t} U(t-r)(\tilde{w}_{\delta,R}(\Gamma+r)+U(\Gamma+r)u_0)\Phi_1 dW^{(1)}_{\Gamma+r}.
\end{aligned}
$$
Then, if we divide $[0,T_0]$ into finite numbers of small intervals 
$$[0,T], [T,2T],... ,[T_0-T, T_0], $$
and using a fixed point argument part by part like \eqref{2_estimate_local}, then we can prove 
$$
w_{\delta,R}\in L_{\omega}^{2m}\left(Z^{T_0}_{b,0,q}\cap C([0,T_0];L_x^2\cap \dot{H}_{x}^{-q})\right), \ \forall R\in(0,+\infty),
$$
and $ u_{\delta,R}\in L_{\omega}^{2m}\left(Y^{T_0}_{b,0}\cap C([0,T_0];L_x^2)\right), \ \forall R\in(0,+\infty)$.
\end{proof}

For the fixed $T_0>0$, let us define the stopping times 
$$
\kappa_{\delta,R}(\omega)=\inf\{t\geq 0: \|u_{\delta,R}\|_{Z^T_{b,0,q}}\geq R\}\wedge T_0, \ \kappa_{\delta}^*=\lim_{R\uparrow \infty}\kappa_{\delta,R}(\omega),\  \mathbb{P}\text{-a.s.}
$$

Since $w_{\delta}=w_{\delta,R}  $ in $[0, \kappa_{\delta,R}\wedge \sigma_{\delta,R}]$, a necessary condition of  
$$w_\delta\in  L_{\omega}^{2m}\left(Z^{T_0}_{b,0,q}\cap C([0,T_0];L_x^2\cap \dot{H}_{x}^{-q})\right)$$
is $ \kappa_\delta^*=T_0$ almost surely.  
To prove this fact, we will use $L_x^2$ conservation law of $u_\delta$. Correspondingly, we also need to smooth the noise kernel and the initial values to overcome the regularity problems in using the conservation law.

\begin{prop}\label{prop:well-posedness_with_moments}
	Under the conditions of Lemma \ref{lem:local_well-posedness_moment_estimate}, we have
		 \begin{equation}\label{bourgain_norm_prior_esti_prop_1}
		\begin{aligned}
			\mathbb{E}\left\| w_\delta\right\|^{2m}_{ C([0,T_0];L_x^2)}
			<C(T_0,m,\|u_0\|_{L_x^2},\|k_1\|_{L_x^2}),
		\end{aligned}
	\end{equation}
	and 
	 \begin{equation}\label{bourgain_norm_prior_esti_prop_2}
		\begin{aligned}
			&\mathbb{E}\left\| w_\delta\right\|^{2m}_{Z^{T_0}_{b,0,q}\cap C([0,T_0];\dot{H}_{x}^{-q})}\\
			<&C(\|u_0\|_{L_x^2},\|v_0\|_{L_x^2},\|k_1\|_{L_x^2\cap L_x^{2/(1+4q)}},\|k_2\|_{L_x^\infty L_y^2},a,b,c, m,T_0,\delta),
		\end{aligned}
	\end{equation}
for any $m\in\mathbb{N}^+$.

What's more, if $\|k_2(x,y)\|_{L_x^\infty L_y^2} $ is less than a constant $C^\dagger(a,c)>0$, then we have a uniform estimate
\begin{equation}\label{uniform_moment_estimate_delta}
 \begin{aligned}
	&\mathbb{E}\left\| w_\delta\right\|^{4}_{Z^{T_0}_{b,0,q}\cap C([0,T_0];L_x^2\cap \dot{H}_{x}^{-q})}\\
	<&C(\|u_0\|_{L_x^2},\|v_0\|_{L_x^2},\|k_1\|_{L_x^2\cap L_x^{2/(1+4q)}},\|k_2\|_{L_x^\infty L_y^2},a,b,c, m,T_0),
\end{aligned}
\end{equation}
with respect to $\delta\in(0,1)$.

\end{prop}
\begin{proof}
	The proof of this proposition can be divided into four steps.
	Some details  in this proof will be omitted, since they are similar to the relative proof in \cite{kdvmuti} and \cite{kdv_sharper}. Additionally, we will outline the proof at the start of each step.
	 
	{\bf Step 1(high regularity local well-posedness of the smoothing equations):} In this step, we need to smooth the initial values and the noise kernels. Then, we will prove that these smoothing equations are in a high regularity spaces.
	For any $n\in \mathbb{N}^+$, $R>0$, $i=\{1,2\}$, let us set 
	$$
	\Phi_i^{(n)}f=\int_{\mathbb{R}} P_nk_i(x-y) f(y)dy.
	$$

	According to Lemma \ref{lem:local_well-posedness_moment_estimate}, Theorem 2.8 in \cite{kdvmuti}, Lemma \ref{lem:well-posedness_dissipative} and Remark \ref{rem:explain_local-wellposedness}, for the following equations
	\begin{equation}\label{smooth_equation_2_1}
	\left\{\begin{aligned}
	dw_{\delta}^{(n)}=&-\left(\partial_{xxx}w_{\delta}^{(n)}+\frac{1}{2}\partial_{x}\left(\left(w_{\delta}^{(n)}+U(t)P_nu_0\right)^2\right)\right)dt\\
	&-\left(w_{\delta}^{(n)}+U(t)P_nu_0\right)\left(v^{(n)}_{\delta}\right)^2dt+\left(w_{\delta}^{(n)}+U(t)P_nu_0\right)\Phi_1^{(n)}dW_t^{(1)},\\
	dv_\delta^{(n)}=&\delta^{-1}\left(\mathcal{L}v_\delta^{(n)}-v_\delta^{(n)}\right)dt+\delta^{-1/2}v_\delta^{(n)}\Phi_2^{(n)} dW_t^2;\\
	w^{(n)}_{\delta}(0) = &0,\ v^{(n)}_{\delta}(0) =P_nv_0,
	\end{aligned}\right.
	\end{equation}
	we have $w_\delta^{(n)}\in C\left(\left[0,\kappa_{\delta,R}^{(n)}\right];H_x^2\right)$ almost surely, where 
	$$
	\kappa_{\delta,R}^{(n)}=\inf\left\{t\geq 0: \left\|w_{\delta}^{(n)} \right\|_{Z^t_{b,0,q}}\geq R \right\}\wedge T_0,\  \mathbb{P}\text{-a.s.}
	$$
	
	{\bf Step 2(proving the $L_\omega^{2m}\left(C\left(\left[0,\kappa_{\delta,R}^{(n)}\right];L_x^2\right)\right)$ prior estimate of $w_\delta^{(n)}$):}
	Since $(w_\delta^{(n)},v_{\delta}^{(n)})\in C\left(\left[0,\kappa_{\delta,R}^{(n)}\right];H_x^2\times H_x^2\right)$, similar like the step 1 of Lemma \ref{lem:priori esimate}, we can prove
	\begin{equation}\label{priori_L_2_2}
	\begin{aligned}
		\mathbb{E}\sup_{t\in[0,\kappa^{(n)}_{\delta,R}]}\left(w_{\delta}^{(n)},w_{\delta}^{(n)}\right)^{m}\leq& C(T_0,m,\|u_0\|_{L_x^2},\|k_1\|_{L_x^2}),
	\end{aligned}
	\end{equation}
	for any $m\in\mathbb{N}^+$.
	
	{\bf Step 3(proving the $L_\omega^{2m}Z^{\kappa_{\delta,R}^{(n)}}_{b,0,q}$ prior estimate of $w_\delta^{(n)}$):}
	In this step, from a path-wise view in in \cite{kdvmuti}, we are going to get the $L_\omega^{2m}Z^{\kappa_{\delta,R}^{(n)}}_{b,0,q}$ prior estimate from $L_\omega^{2m}\left(C\left(\left[0,\kappa_{\delta,R}^{(n)}\right];L_x^2\right)\right)$. Relative reasons why Theorem \ref{thm:global_well-posedness} can not get the well-posedness with moment estimates have been explained in the introduction of \cite{kdvas}.

	According to the proof of Lemma \ref{lem:local_well-posedness_moment_estimate}, under the corresponding restriction of constants, for any $R>0$, $t(\omega)\in \left[0,\kappa^{(n)}_{\delta,R}\wedge 1\right]$ and $m\in\mathbb{N}^+$, we have 
	$$
	\begin{aligned}
	&\left\|w_{\delta}^{(n)}\right\|^{2m}_{Z^{t(\omega)}_{b,0,q}\cap C([0,t(\omega)];L_x^2\cap \dot{H}_{x}^{-q})}\\
	\lesssim & \left((t(\omega))^{1-(a+c)}+(t(\omega))^{1/4-a/2}\right)^{2m}\left(\left\|w_\delta^{(n)}\right\|^{4m}_{Z^{t(\omega)}_{b,0,q}}+\|u_0\|_{L_x^2}^{4m} \right) \\
	&+\left((t(\omega))^{1-(a+c)}+(t(\omega))^{1/4-a/2}\right)^{2m}\left(\left\|w_\delta^{(n)}\right\|^{2m}_{Z^{t(\omega)}_{b,0,q}}+\|u_0\|_{L_x^2}^{2m} \right)\left\| v_{\delta}^{(n)}\right\|^{4m}_{L_{t(\omega)}^4L_x^2}\\
	&+\left\|\int_{0}^{t} U(t-r)\left(\left(w_{\delta}^{(n)}+U(r)P_nu_0\right)\Phi_1^{(n)} dW^{(1)}_r\right)\right\|^{2m}_{Z^{t(\omega)}_{b,0,q}\cap C([0,t(\omega)];L_x^2\cap \dot{H}_{x}^{-q})}.
	\end{aligned}
	$$
	
	Just like Lemma 3.2 in \cite{kdvmuti}, let us set
	$$
	\begin{aligned}
	K_1(\omega)=&\left\|\int_{0}^{t} U(t-r)\left(\left(w_{\delta}^{(n)}+U(r)P_nu_0\right)\Phi_1^{(n)} dW^{(1)}_r\right)\right\|^{2m}_{Z^{\kappa^{(n)}_{\delta,R}}_{b,0,q}\cap C\left([0,\kappa^{(n)}_{\delta,R}];L_x^2\cap \dot{H}_{x}^{-q}\right)}\\
	&+ \left\|w_{\delta}^{(n)}\right\|^{6m}_{C\left(\left[0,\kappa_{\delta,R}^{(n)}\right];L_x^2\right)}+\|u_0\|^{6m}_{L_x^2}+\left\|v_{\delta}^{(n)}\right\|^{6m}_{L_{T_0}^4L_x^2}
	\end{aligned}
	$$ 
	Thus, if we choose $t(\omega)$ sufficiently small such that 
	\begin{equation}\label{interval_number_1}
		1/16<\left(t(\omega))^{1-(a+c)}+(t(\omega))^{1/4-a/2}\right)K_1(\omega)<1/8
	\end{equation}
	and 
	\begin{equation}\label{interval_number_2}
1/16<\left(t(\omega))^{1-(a+c)}+(t(\omega))^{1/4-a/2}\right)\left\|v_\delta^{(n)}\right\|^2_{L_{T_0}^4L_x^2}<1/8,
	\end{equation}
	then we have
	$$
	\begin{aligned}
		\left\|w_{\delta}^{(n)}\right\|^{2m}_{Z^{t(\omega)}_{b,0,q}\cap C([0,t(\omega)];L_x^2\cap \dot{H}_{x}^{-q})}
		\lesssim 2 K_1(\omega).
	\end{aligned}
	$$
	
	Let us  divide  $\left[0,\kappa_{\delta,R}^{(n)}\right]$ into finite numbers of small intervals 
	$$[0,t(\omega)],\ [t(\omega),2t(\omega)],\ ...,\ \left[(l-1)t(\omega),\kappa^{(n)}_{\delta,R}\right]$$
	and set 
	\begin{equation}\label{aux_linear}	w_{\delta}^{(n,k)}(t)=u_{\delta}^{(n)}-U(t-kt(\omega))u_\delta^{(n)}(kt(\omega)),\ \forall t\in[kt(\omega),(k+1)t(\omega)],
	\end{equation}
	for any $k\in\{0,1,2,...,l-1\}$.
	
	Similar like above, we have 
	$$
	\left\|w_{\delta}^{(n,k)}\right\|^{2m}_{Z^{[kt(\omega),(k+1)t(\omega)]}_{b,0,q}\cap C([kt(\omega),(k+1)t(\omega)];L_x^2\cap \dot{H}_{x}^{-q})}\\
	\lesssim 2K_1(\omega).
	$$
	
	What's more, for any $m\in\mathbb{N}^+$, we have
	$$
	\begin{aligned}
	&\left\|w_{\delta}^{(n)}\right\|^{2m}_{C([0,\kappa_{\delta,R}^{(n)}];\dot{H}_x^{-q})}\\
	\lesssim&\sum_{k=0}^{l-1}\left((t(\omega))^{1-(a+c)}\right)^{2m}\left\|w_\delta^{(n,k)}\right\|^{4m}_{Z^{[kt(\omega),(k+1)t(\omega)]}_{b,0,q}} \\
	&+\sum_{k=0}^{l-1}\left((t(\omega))^{1-(a+c)}\right)^{2m}\left\|w_\delta^{(n,k)}\right\|^{2m}_{Z^{[kt(\omega),(k+1)t(\omega)]}_{b,0,q}} \left\| v_{\delta}^{(n)}\right\|^{4m}_{L_{T_0}^4L_x^2}\\
	&+\left\|\int_{0}^{t} U(t-r)\left(u_{\delta}^{(n)}\Phi_1^{(n)} dW^{(1)}_r\right)\right\|^{2m}_{C([0,\kappa_{\delta,R}^{(n)}]; \dot{H}_{x}^{-q})}+lK_1(\omega)\\
	\lesssim&lK_1(\omega)+\left\|\int_{0}^{t} U(t-r)\left(u_{\delta}^{(n)}\Phi_1^{(n)} dW^{(1)}_r\right)\right\|^{2m}_{C([0,\kappa_{\delta,R}^{(n)}]; \dot{H}_{x}^{-q})}\\
	:=&K_2.
	\end{aligned}
	$$
	
	Furthermore, according to \eqref{aux_linear}, we have 
	$$
    w_{\delta}^{(n)}=w_{\delta}^{(n,k)}(t)+U(t-kt(\omega))u_\delta^{(n)}(kt(\omega))-U(t)P_nu_0,\ \forall t\in[kt(\omega),(k+1)t(\omega)].
	$$ 
	
	Thus, it can be deduce from the definitions of $K_1(\omega),K_2(\omega)$ that 
	$$
	\begin{aligned}
    &\left\| w_\delta^{(n)}\right\|^{2m}_{Z^{\kappa_{\delta}^{(n)}}_{b,0,q}\cap C([0,\kappa_{\delta,R}^{(n)}];L_x^2\cap \dot{H}_{x}^{-q})}\\
    \lesssim& \sum_{k=0}^{l-1} \left\| w_\delta^{(n)}\right\|^{2m}_{Z^{[kt(\omega),(k+1)t(\omega)]}_{b,0,q}\cap C([kt(\omega),(k+1)t(\omega)];L_x^2\cap \dot{H}_{x}^{-q})}\\
    \lesssim&\sum_{k=0}^{l-1}K_1(\omega)+K_2(\omega)+\|u_0\|_{L_x^2}^{2m}+\left\|u_\delta^{(n)}\right\|^{2m}_{C([0,\kappa_{\delta,R}^{(n)}];L_x^2)},
	\end{aligned}
	$$
	for any $m\in\mathbb{N}^+$. 
	
	Therefore, we have 
	$$
	\begin{aligned}
	&\mathbb{E}\left\| w_\delta^{(n)}\right\|^{2m}_{Z^{\kappa_{\delta,R}^{(n)}}_{b,0,q}\cap C([0,\kappa_{\delta,R}^{(n)}];L_x^2\cap \dot{H}_{x}^{-q})}\\
	\lesssim&\mathbb{E}(lK_1(\omega))+\mathbb{E}(l^2K_1(\omega))+\mathbb{E}\left(l\left\|u_\delta^{(n)}\right\|^{2m}_{C([0,\kappa_{\delta,R}^{(n)}];L_x^2)}\right)\\
	&+\mathbb{E}\left(l\left\|\int_{0}^{t} U(t-r)\left(u_{\delta}^{(n)}\Phi_1^{(n)} dW^{(1)}_r\right)\right\|^{2m}_{C([0,\kappa_{\delta,R}^{(n)}]; \dot{H}_{x}^{-q})}\right).
	\end{aligned}
	$$
	According to \eqref{interval_number_1} and \eqref{interval_number_2}, we have 
	$$
	\begin{aligned}
	l\lesssim &T_0\cdot (K_1(\omega))^{1/(1-(a+c))} \vee T_0\cdot (K_1(\omega))^{4/(1-2a)} \\
	&\vee T_0\cdot \left\|v_\delta^{(n)}\right\|_{L_{T_0}^4L_x^2}^{2/(1-(a+c))}  \vee T_0\cdot \left\|v_\delta^{(n)}\right\|_{L_{T_0}^4L_x^2}^{8/(1-2a)}.
	\end{aligned}
	$$
	
	Thus, based on Lemma \ref{lem:stochastic_convolution_estimate_2}, \eqref{priori_L_2_2} and Lemma \ref{lem:well-posedness_dissipative}, by the H\"{o}lder inequality, we have
	\begin{equation}\label{bourgain_norm_prior_esti}
	\begin{aligned}
		&\mathbb{E}\left\| w_\delta^{(n)}\right\|^{2m}_{Z^{\kappa_{\delta,R}^{(n)}}_{b,0,q}\cap C([0,\kappa_{\delta,R}^{(n)}];L_x^2\cap \dot{H}_{x}^{-q})}\\
		\lesssim& \mathbb{E}\left\|u_{\delta}^{(n)}\right\|_{C([0,\kappa_{\delta,R}^{(n)}];L_x^2)}^{C(a,c,m)}+\mathbb{E}\left\|v_{\delta}^{(n)}\right\|^{C(a,c,m)}_{L_{T_0}^{4}L_x^2}\\
		&+\mathbb{E}\left\|\int_{0}^{t} U(t-r)\left(\left(w_{\delta}^{(n)}+U(r)P_nu_0\right)\Phi_1^{(n)} dW^{(1)}_r\right)\right\|^{C(a,c,m)}_{Z^{\kappa^{(n)}_{\delta,R}}_{b,0,q}\cap C\left([0,\kappa^{(n)}_{\delta,R}];L_x^2\cap \dot{H}_{x}^{-q}\right)}\\
		\lesssim& C(\|k_1\|_{L_x^2\cap L_x^{2/(1+4q)}},a,b,c, m,T_0)\cdot\mathbb{E}\left\|u_{\delta}^{(n)}\right\|_{C([0,\kappa_{\delta,R}^{(n)}];L_x^2)}^{C(a,c,m)}+\mathbb{E}\left\|v_{\delta}^{(n)}\right\|^{C(a,c,m)}_{L_{T_0}^{4}L_x^2}\\
		\lesssim&  C(\|u_0\|_{L_x^2},\|v_0\|_{L_x^2},\|k_1\|_{L_x^2\cap L_x^{2/(1+4q)}},\|k_2\|_{L_x^\infty L_y^2},a,b,c, m,T_0,\delta),
	\end{aligned}
	\end{equation}
 where $C(a,c,m)$ is a constant only depends on $a,c,m$.
 
 Hence, by the monotone convergence theorem,  \eqref{bourgain_norm_prior_esti} implies
 $$
 \begin{aligned}
  &\lim_{R\uparrow +\infty}\mathbb{E}\left\| w_\delta^{(n)}\right\|^{2m}_{Z^{\kappa_{\delta,R}^{(n)}}_{b,0,q}\cap C([0,\kappa_{\delta,R}^{(n)}];L_x^2\cap \dot{H}_{x}^{-q})}\\
  \lesssim& C(\|u_0\|_{L_x^2},\|v_0\|_{L_x^2},\|k_1\|_{L_x^2\cap L_x^{2/(1+4q)}},\|k_2\|_{L_x^\infty L_y^2},a,b,c, m,T_0,\delta)<+\infty
 \end{aligned}
 $$
 and 
 $\lim_{R\uparrow +\infty} \kappa_{\delta,R}^{(n)}=T_0$ almost surely.
 
 Namely, we have 
 \begin{equation}\label{bourgain_norm_prior_esti_final}
 \begin{aligned}
  &\mathbb{E}\left\| w_\delta^{(n)}\right\|^{2m}_{Z^{T_0}_{b,0,q}\cap C([0,T_0];L_x^2\cap \dot{H}_{x}^{-q})}\\
  <&C(\|u_0\|_{L_x^2},\|v_0\|_{L_x^2},\|k_1\|_{L_x^2\cap L_x^{2/(1+4q)}},\|k_2\|_{L_x^\infty L_y^2},a,b,c, m,T_0,\delta),
  \end{aligned}
 \end{equation}
for any $n,m\in\mathbb{N}^+$.
 
What's more, according to Lemma  \ref{lem:well-posedness_dissipative}, there exists a constant $C^{\dagger}$ only depending on $C(a,c,2)$ in \eqref{bourgain_norm_prior_esti}, such that once $ \|k_2\|_{L_x^\infty L_y^2}<C^{\dagger} $, then we have $ \mathbb{E}\left\|v_{\delta}^{(n)}\right\|^{C(a,c,2)}_{L_{T_0}^{4}L_x^2}\lesssim  C(T_0,a,c,\|v_0\|_{L_x^2})$, which finally implies 
$$
 \begin{aligned}
	&\mathbb{E}\left\| w_\delta^{(n)}\right\|^{4}_{Z^{T_0}_{b,0,q}\cap C([0,T_0];L_x^2\cap \dot{H}_{x}^{-q})}\\
	<&C(\|u_0\|_{L_x^2},\|v_0\|_{L_x^2},\|k_1\|_{L_x^2\cap L_x^{2/(1+4q)}},\|k_2\|_{L_x^\infty L_y^2},a,b,c, m,T_0),
\end{aligned}
$$
uniformly with respect to $\delta\in(0,1)$.

{\bf Step 4(proving the $L_\omega^{2m}\left(Z^{T_0}_{b,0,q}\cap C([0,T_0];L_x^2\cap \dot{H}_{x}^{-q})\right)$ prior estimate of $w_\delta$ by letting $n\uparrow+\infty$):}
In this step, we will prove \eqref{bourgain_norm_prior_esti_prop_1} and \eqref{bourgain_norm_prior_esti_prop_2} through a approximation method. A lot of proof will be omitted, since it is similar to  step 2 and step 3 in Lemma \ref{lem:priori esimate} or Theorem 1.2 in \cite{kdv_sharper}.

Let us define $\tilde{w}_\delta^{(n)}= w_\delta-w_\delta^{(n)}$ and
$$
\tilde{\kappa}_{\delta,R}^{(n)}=\kappa_{\delta,R}^{(n)}\wedge\kappa_{\delta,R}\wedge\inf\{t\geq 0:\|v_\delta\|_{L_t^4L_x^2}\geq R\} \wedge \inf\{t\geq 0:\|v^{(n)}_\delta\|_{L_t^4L_x^2}\geq R\}
$$
almost surely.

We will study the equation 
\begin{equation*}
\begin{aligned}
&\tilde{w}_\delta^{(n)}\left(\iota\right)\\
=&\int_{0}^{\iota}U(\iota-r)\partial_{x}\left( (\tilde{w}_\delta^{(n)}+U(r)P_{>n}u_0)(w_\delta+w_\delta^{(n)}+U(r)P_{n}u_0+U(r)u_0)\right)dr\\
&+\int_{0}^{\iota}U(\iota-r)\left((\tilde{w}_\delta^{(n)}+U(r)P_{>n}u_0)\left(v_{\delta}^{(n)}\right)^2\right)dr\\
&+\int_{0}^{\iota}U(\iota-r)\left((w_\delta+U(r)u_0)\left(v_{\delta}^{(n)}+v_\delta\right)\left(v_\delta-v_\delta^{(n)}\right)\right)dr\\
&+\int_{0}^{\iota} U(\iota-r)\left( (\tilde{w}_{\delta}+U(r)P_{>n}u_0 )\Phi_1dW^{1}_r\right)\\
&+\int_{0}^{\iota} U(\iota-r)\left( (w_{\delta}^{(n)}+U(r)P_{n}u_0 )\left(\Phi_1\right)_{>n}dW^{1}_r\right),
\end{aligned}
\end{equation*}
for any stopping time satisfying $ 0\leq \iota\leq\tilde{\kappa}_{\delta,R}^{(n)} $ almost surely.

Similar like before, for any $T\in[0,T_0]$, $m\in\mathbb{N}^+$, we have 
\begin{equation}\label{difference_n}
	\begin{aligned}
		&\left\|\tilde{w}_\delta^{(n)}\right\|_{Z^{T\wedge\tilde{\kappa}_{\delta,R}^{(n)} }_{b,0,q}\cap C([0,T\wedge\tilde{\kappa}_{\delta,R}^{(n)}];L_x^2\cap \dot{H}_{x}^{-q})}\\
		\lesssim& \left(T^{1-(a+c)}+T^{1/4-a/2}\right)R\left\|\tilde{w}_\delta^{(n)}\right\|_{Z^{T\wedge\tilde{\kappa}_{\delta,R}^{(n)} }_{b,0,q}}\\
		&+ \left(T^{1-(a+c)}+T^{1/4-a/2}\right)\left(\|u_0\|_{L_x^2}\|P_{>n}u_0\|_{L_x^2}+R\|P_{>n}u_0\|_{L_x^2}\right)\\
		&+\left(T^{1-(a+c)}+T^{1/4-a/2}\right)R^2\left\|\tilde{w}_\delta^{(n)}\right\|_{Z^{T\wedge\tilde{\kappa}_{\delta,R}^{(n)} }_{b,0,q}}\\
    	&+\left(T^{1-(a+c)}+T^{1/4-a/2}\right)R^2\left\|P_{>n}u_0\right\|_{L_x^2} \\
    	&+ \left(T^{1-(a+c)}+T^{1/4-a/2}\right)\left(R^2+\|u_0\|^2_{L_x^2}\right)\|v_{\delta}^{(n)}-v_\delta\|_{L_{T_0}^4L_x^2}\\
		&+\left\|\int_{0}^{\iota} U(\iota-r)\left( (\tilde{w}_{\delta}+U(r)P_{>n}u_0 )\Phi_1dW^{1}_r\right)\right\|_{Z^{T\wedge\tilde{\kappa}_{\delta,R}^{(n)} }_{b,0,q}\cap C([0,T\wedge\tilde{\kappa}_{\delta,R}^{(n)}];L_x^2\cap \dot{H}_{x}^{-q})}\\
		&+\left\|\int_{0}^{\iota} U(\iota-r)\left( (w_{\delta}^{(n)}+U(r)P_{n}u_0 )\left(\Phi_1\right)_{>n}dW^{1}_r\right)\right\|_{Z^{T\wedge\tilde{\kappa}_{\delta,R}^{(n)} }_{b,0,q}\cap C([0,T\wedge\tilde{\kappa}_{\delta,R}^{(n)}];L_x^2\cap \dot{H}_{x}^{-q})}.
	\end{aligned}
\end{equation}

Thus, just like Lemma \ref{lem:priori esimate}, if we divide $[0,\tilde{\kappa}_{\delta,R}^{(n)}]$ into finite numbers of intervals smaller than a $ T^*(a,c,R,\|u_0\|_{L_x^2},\|k_1\|_{L_x^2\cap L_x^{2/(1+4q)}})$  and let $n\uparrow+\infty$, we can prove \eqref{bourgain_norm_prior_esti_prop_1},  \eqref{bourgain_norm_prior_esti_prop_2} and \eqref{uniform_moment_estimate_delta}.

\end{proof}

In the proof of Theorem \ref{thm:result_with_classical_workspace}, owing to the multilinear estimates of the nonlinear terms,  we need to divide $[0,T_0]$ into two parts to get an estimate of $\mathbb{E}\|u-u_\delta\|^2_{C([0,T_0];L_x^2)}$. Thus, we introduce the following stopping times
$$
\iota_{\delta,R}:=\inf\left\{t\geq0:\| w_\delta\|_{Z^t_{b,0,q}\cap C\left([0,t];L_x^2\cap L_x^{2/(1+4q)}\right)}\geq R\right\},\ \mathbb{P}\text{-a.s.}
$$
for any $ R>0$, $\delta\in[0,1)$. Here, $w_0(t):=u(t)-U(t)u_0$ is the solution of \eqref{average_eq} subtracting the linear part.

To get a uniform tail estimate of $\mathbb{E}\|u_\delta\|^2_{Z^{[T_0\wedge \iota_{\delta,R},T_0]}_{b,0,q}} $, we have the following corollary of Proposition \ref{prop:well-posedness_with_moments}.

\begin{coro}\label{coro:uniform_tail_estimate}
Under the conditions of Proposition \ref{prop:well-posedness_with_moments}, we have  
\begin{equation}\label{uniform_tail_esti}
	\begin{aligned}
	&\mathbb{E}\|w_\delta\|^2_{Z^{[T_0\wedge \iota_{\delta,R},T_0]}_{b,0,q}\cap C([T_0\wedge \iota_{\delta,R},T_0];L_x^2\cap \dot{H}_{x}^{-q})}\\
	\leq& \frac{	\mathbb{E}\|w_\delta\|^4_{Z^{T_0}_{b,0,q}\cap C([0,T_0];L_x^2\cap \dot{H}_{x}^{-q})}}{R^2}\\
	\leq& C(\|u_0\|_{L_x^2},\|v_0\|_{L_x^2},\|k_1\|_{L_x^2\cap L_x^{2/(1+4q)}},\|k_2\|_{L_x^\infty L_y^2},a,b,c,T_0)/R^2, 
	\end{aligned}
\end{equation}
for any $\delta\in[0,1)$.
\end{coro}
\begin{proof}
	
	Note that to prove \eqref{uniform_tail_esti}, it is sufficient to prove 
	\begin{equation}\label{trick}
	R^2\|w_\delta\|^2_{Z^{[T_0\wedge \iota_{\delta,R},T_0]}_{b,0,q}\cap C([T_0\wedge \iota_{\delta,R},T_0];L_x^2\cap \dot{H}_{x}^{-q})}\leq \|w_\delta\|^4_{Z^{T_0}_{b,0,q}\cap C([0,T_0];L_x^2\cap \dot{H}_{x}^{-q})},
	\end{equation}
    for any $  \omega\in\Omega$.
    
    In fact, on one hand, for any $\omega\in \Omega$ such that $ \|w_\delta\|_{Z^{T_0}_{b,0,q}\cap C([0,T_0];L_x^2\cap \dot{H}_{x}^{-q})}\leq R$, we have 
    $$ 
    \begin{aligned}
    &R^2\|w_\delta\|^2_{Z^{[T_0\wedge \iota_{\delta,R},T_0]}_{b,0,q}\cap C([T_0\wedge \iota_{\delta,R},T_0];L_x^2\cap \dot{H}_{x}^{-q})}=0\\
    \leq& \|w_\delta\|^4_{Z^{T_0}_{b,0,q}\cap C([0,T_0];L_x^2\cap \dot{H}_{x}^{-q})}.
    \end{aligned}
    $$ 
    
    On the other hand, for any $\omega\in \Omega$ such that $ \|w_\delta\|_{Z^{T_0}_{b,0,q}\cap C([0,T_0];L_x^2\cap \dot{H}_{x}^{-q}) }> R$, according to the definition of restricted norms, we also have 
    $$ 
    \|w_\delta\|_{Z^{[T_0\wedge \iota_{\delta,R},T_0]}_{b,0,q}\cap C([T_0\wedge \iota_{\delta,R},T_0];L_x^2\cap \dot{H}_{x}^{-q})}\leq  \|w_\delta\|_{Z^{T_0}_{b,0,q}\cap C([0,T_0];L_x^2\cap \dot{H}_{x}^{-q})}.
    $$ 
    
    Combining \eqref{uniform_moment_estimate_delta} with the fact  
    \begin{equation}\label{well-posedness_w}
    	w(t):=u(t)-U(t)u_0\in L^{2m}\left(Z^{T_0}_{b,0,q} \cap C\left([0,T_0];L_x^2\cap H_x^{-q}\right)\right),
    \end{equation}
    for any $m\in\mathbb{N}^+$ in \cite{kdvmuti} or \cite{kdv_sharper},  
    we can illustrate  \eqref{trick} is true.

\end{proof}

Now, we are in the position to prove Theorem \ref{thm:result_with_classical_workspace}.
\begin{proof}[{\bf Proof of Theorem \ref{thm:result_with_classical_workspace}:}]

Let us define $\tilde{w}_\delta=w_0-w_\delta$,
$$
\tilde{\iota}_{\delta,R}= \iota_{0,R}\wedge \iota_{\delta,R}, \ \mathbb{P}\text{-a.s.}
$$

According to Corollary \ref{coro:uniform_tail_estimate}, for any $\varepsilon >0$, there exists a  $R>\|u_0\|_{L_x^2}\vee1$  such that 
\begin{equation*}
\begin{aligned}
&\mathbb{E}\|w_\delta\|^2_{Z^{[T_0\wedge \iota_{\delta,R},T_0]}_{b,0,q}}+\mathbb{E}\|w_0\|^2_{Z^{[T_0\wedge\iota_{0,R},T_0]}_{b,0,q}}\\
\leq &\frac{4	\mathbb{E}\|w_\delta\|^4_{Z^{T_0}_{b,0,q}\cap C([0,T_0];L_x^2\cap \dot{H}_{x}^{-q})}}{R^2}+\frac{4	\mathbb{E}\|w_0\|^4_{Z^{T_0}_{b,0,q}\cap C([0,T_0];L_x^2\cap \dot{H}_{x}^{-q})}}{R^2} < \varepsilon,
\end{aligned}
\end{equation*}
for any $ \delta\in(0,1)$. Then, we have 
\begin{equation}\label{tail_moment_estimate}
	\begin{aligned}
		&\mathbb{E}\|w_\delta\|^2_{Z^{[T_0\wedge \iota_{\delta,\sqrt{l}R/2},T_0]}_{b,0,q}}+\mathbb{E}\|w_0\|^2_{Z^{[T_0\wedge\iota_{\sqrt{l}R/2},T_0]}_{b,0,q}}\\
		\leq &\frac{4	\mathbb{E}\|w_\delta\|^4_{Z^{T_0}_{b,0,q}\cap C([0,T_0];L_x^2\cap \dot{H}_{x}^{-q})}}{lR^2}+\frac{4	\mathbb{E}\|w\|^4_{Z^{T_0}_{b,0,q}\cap C([0,T_0];L_x^2\cap \dot{H}_{x}^{-q})}}{lR^2} < \varepsilon/l,
	\end{aligned}
\end{equation}
for any $\delta\in(0,1)$ and $l>0$.

We will study the following equation 
$$
\begin{aligned}
\tilde{w}_\delta(\iota)
=&\int_{0}^{\iota}U(\iota-r)\partial_{x}\left( \tilde{w}_\delta(w_\delta+w_0+2U(r)u_0)\right)dr\\
&+\int_{0}^{\iota}U(\iota-r)\left( (w_\delta+U(r)u_0)v_\delta^2\right)dr\\
&+\int_{0}^{\iota} U(\iota-r)\left( \tilde{w}_{\delta} \Phi_1dW^{1}_r\right),
\end{aligned}
$$
for any stopping time $\iota$ satisfying $0\leq \iota \leq \tilde{\iota}_{\delta,\sqrt{l}R}$ almost surely.

Thus, according to the definition of $\tilde{\iota}_{\delta,\sqrt{l}R} $, for any fixed $R>\|u_0\|_{L_x^2}$ and $l>1$, we have 
$$
\begin{aligned}
&\mathbb{E}\|\tilde{w}_\delta\|^2_{Z^{t\wedge \tilde{\iota}_{\delta,\sqrt{l}R}}_{b,0,q}\cap C([0,t\wedge \tilde{\iota}_{\delta,\sqrt{l}R}];L_x^2\cap \dot{H}_{x}^{-q}) }
\\
\lesssim& \left(t^{2-2(a+c)}+t^{1/2-a}\right) \left(\|u_0\|_{L_x^2}^2+lR^2\right)\mathbb{E}\|\tilde{w}_\delta\|^2_{Z^{t\wedge \tilde{\iota}_{\delta,\sqrt{l}R}}_{b,0,q} }\\
&+\left(t^{2-2(a+c)}+t^{1/2-a}\right)\left(\|u_0\|_{L_x^2}^2+lR^2\right)\mathbb{E}\|v_\delta\|^4_{L_{T_0}^4L_x^2}\\
&+t\cdot(\|k_1\|^{2}_{L_x^{2/(1+4q)}}+\|k_1\|_{L_x^2}^{2})\mathbb{E}\|\tilde{w}_\delta\|^{2}_{C([0,t\wedge \tilde{\iota}_{\delta,\sqrt{l}R}];L_x^2)}\\
\lesssim& \left(t^{2-2(a+c)}+t^{1/2-a}\right) lR^2\left(\mathbb{E}\|\tilde{w}_\delta\|^2_{Z^{t\wedge \tilde{\iota}_{\delta,\sqrt{l}R}}_{b,0,q} }+\mathbb{E}\|v_\delta\|^4_{L_{T_0}^4L_x^2}\right)\\
&+t\cdot(\|k_1\|^{2}_{L_x^{2/(1+4q)}}+\|k_1\|_{L_x^2}^{2})\mathbb{E}\|\tilde{w}_\delta\|^{2}_{C([0,t\wedge \tilde{\iota}_{\delta,\sqrt{l}R}];L_x^2)},
\end{aligned}
$$
for any $t\in[0,T_0 \wedge 1]$.

Choosing $T^*=\left(l^{1/(2a-1)}\wedge l^{1/(2(a+c)-2)}\right)\cdot T^\dagger(R,\|k_1\|_{L_x^2\cap L_x^{2/(1+4q)}})$ sufficiently small, we have   
$$
\begin{aligned}
&\mathbb{E}\|\tilde{w}_\delta\|^2_{Z^{T^*\wedge \tilde{\iota}_{\delta,\sqrt{l}R}}_{b,0,q}\cap C([0,T^*\wedge \tilde{\iota}_{\delta,\sqrt{l}R}];L_x^2\cap \dot{H}_{x}^{-q}) }\\
\leq& \frac{1}{4}\left(\mathbb{E}\|\tilde{w}_\delta\|^2_{Z^{T^*\wedge \tilde{\iota}_{\delta,\sqrt{l}R}}_{b,0,q}\cap C([0,T^*\wedge \tilde{\iota}_{\delta,\sqrt{l}R}];L_x^2\cap \dot{H}_{x}^{-q}) }+\mathbb{E}\|v_\delta\|^4_{L_{T_0}^4L_x^2}\right),
\end{aligned}
$$
which means $$\mathbb{E}\|\tilde{w}_\delta\|^2_{Z^{T^*\wedge \tilde{\iota}_{\delta,\sqrt{l}R}}_{b,0,q}\cap C([0,T^*\wedge \tilde{\iota}_{\delta,\sqrt{l}R}];L_x^2\cap \dot{H}_{x}^{-q}) }\leq \delta,$$
by Lemma \ref{lem:well-posedness_dissipative}.  
Here,  $T^\dagger$ is a constant only depending on  $R$ and  $\|k_1\|_{L_x^2\cap L_x^{2/(1+4q)}}$.

Dividing $ [0,T_0\wedge\tilde{\iota}_{\delta,\sqrt{l}R} ]$ into finite numbers $K$ ($K$ is directly proportional to $ l^{1/(1-2a)}\vee l^{1/(2-2(a+c))}$) of small intervals 
$$
[0,T^*\wedge \tilde{\iota}_{\delta,\sqrt{l}R}],...,\ [(k-1)T^*\wedge \tilde{\iota}_{\delta,\sqrt{l}R},kT^*\wedge \tilde{\iota}_{\delta,\sqrt{l}R} ],...,\ [ (T_0-T^*)\wedge \tilde{\iota}_{\delta,\sqrt{l}R},T_0\wedge \tilde{\iota}_{\delta,\sqrt{l}R}],
$$
and setting
$$
w_{\delta}^{(k)}(t)=u_{\delta}\left(t+\left[(k-1)T^*\wedge \tilde{\iota}_{\delta,\sqrt{l}R}\right]\right)-U\left(t\right)u_\delta\left((k-1)T^*\wedge \tilde{\iota}_{\delta,\sqrt{l}R}\right),
$$
$$
w^{(k)}_0(t)=u\left(t+\left[(k-1)T^*\wedge \tilde{\iota}_{\delta,\sqrt{l}R}\right]\right)-U\left(t\right)u\left((k-1)T^*\wedge \tilde{\iota}_{\delta,\sqrt{l}R}\right),
$$
for any $t\in\left[0,\left(\tilde{\iota}_{\delta,\sqrt{l}R}\wedge k T^*\right)-\left(\tilde{\iota}_{\delta,\sqrt{l}R}\wedge (k-1) T^*\right)\right]$.

It is clear that 
$$
\begin{aligned}
&\tilde{w}_{\delta}(t)\\
&=
w^{(k)}_0\left(t-t-\left[(k-1)T^*\wedge \tilde{\iota}_{\delta,\sqrt{l}R}\right]\right)-w_{\delta}^{(k)}\left(t-t-\left[(k-1)T^*\wedge \tilde{\iota}_{\delta,\sqrt{l}R}\right]\right)\\
&+U\left(t-\left[(k-1)T^*\wedge \tilde{\iota}_{\delta,\sqrt{l}R}\right]\right)\left(u\left((k-1)T^*\wedge \tilde{\iota}_{\delta,\sqrt{l}R}\right)-u_\delta\left((k-1)T^*\wedge \tilde{\iota}_{\delta,\sqrt{l}R}\right)\right) ,
\end{aligned}
$$
for any $ t\in \left[(k-1)T^*\wedge \tilde{\iota}_{\delta,\sqrt{l}R},kT^*\wedge \tilde{\iota}_{\delta,\sqrt{l}R} \right]$ and $k\in\{1,2,...,K\}$.

What's more, we have
$$
\begin{aligned}
&w^{(k)}_0(\iota)-w_{\delta}^{(k)}(\iota)\\
=&\int_{0}^{\iota}U(\iota-r)\partial_{x}\Bigg( \tilde{w}_\delta\big(r+[(k-1)T^*\wedge \tilde{\iota}_{\delta,\sqrt{l}R}]\big)\\
&\cdot\Big(w_\delta(r+[(k-1)T^*\wedge \tilde{\iota}_{\delta,\sqrt{l}R}]\big)+w_0(r+[(k-1)T^*\wedge \tilde{\iota}_{\delta,\sqrt{l}R}]\big)\\
&+2U((r+[(k-1)T^*\wedge \tilde{\iota}_{\delta,\sqrt{l}R}]))u_0\Big)\Bigg)dr\\
&+\int_{0}^{\iota}U(\iota-r)\Bigg( (w_\delta(r+[(k-1)T^*\wedge \tilde{\iota}_{\delta,\sqrt{l}R}])\\
&+U(r+[(k-1)T^*\wedge \tilde{\iota}_{\delta,\sqrt{l}R}])u_0)v_\delta^2(r+[(k-1)T^*\wedge \tilde{\iota}_{\delta,\sqrt{l}R}])\Bigg)dr\\
&+\int_{0}^{\iota} U(\iota-r)\left( \tilde{w}_{\delta}(r+[(k-1)T^*\wedge \tilde{\iota}_{\delta,\sqrt{l}R}]) \Phi_1dW^{1}_r\right),
\end{aligned}
$$
for any stopping times $\iota$ satisfying $0\leq \iota \leq \left(\tilde{\iota}_{\delta,\sqrt{l}R}\wedge k T^*\right)-\left(\tilde{\iota}_{\delta,\sqrt{l}R}\wedge (k-1) T^*\right)$ almost surely.

Thus, if we set 
$$
\alpha_k=\mathbb{E}\| \tilde{w}_{\delta}\|^2_{Z^{[(k-1)T^*\wedge \tilde{\iota}_{\delta,\sqrt{l}R},kT^*\wedge \tilde{\iota}_{\delta,\sqrt{l}R}]}_{b,0,q}\cap C([(k-1)T^*\wedge \tilde{\iota}_{\delta,\sqrt{l}R},kT^*\wedge \tilde{\iota}_{\delta,\sqrt{l}R}];L_x^2\cap \dot{H}_{x}^{-q}) },
$$ 
then we have $\alpha_1\leq \delta$ and 
$$
\begin{aligned}
\alpha_k
\leq&2\mathbb{E}\left\|\tilde{w}_{\delta}((k-1)T^*\wedge \tilde{\iota}_{\delta,\sqrt{l}R})\right\|^2_{L_x^2\cap \dot{H}_{x}^{-q}}\\
& +2\mathbb{E}\left\|w^{(k)}_{\delta}-w_0^{(k)}\right\|^2_{Z^{[(k-1)T^*\wedge \tilde{\iota}_{\delta,\sqrt{l}R},kT^*\wedge \tilde{\iota}_{\delta,\sqrt{l}R} ]}_{b,0,q}\cap C([(k-1)T^*\wedge \tilde{\iota}_{\delta,\sqrt{l}R},kT^*\wedge \tilde{\iota}_{\delta,\sqrt{l}R} ];L_x^2\cap \dot{H}_{x}^{-q})}\\
\leq&2\alpha_{k-1}+\frac{1}{2}\alpha_k+\frac{1}{2}\delta,
\end{aligned}
$$
which implies
$$
\alpha_k\leq 4\alpha_{k-1}+\delta,\ k\in\{2,3,...,K\}.
$$

Thus, we have 
\begin{equation*}
\begin{aligned}
&\mathbb{E}\|\tilde{w}_\delta\|^2_{Z^{T_0\wedge \tilde{\iota}_{\delta,R}}_{b,0,q}\cap C([0,T_0\wedge \tilde{\iota}_{\delta,R}];L_x^2\cap \dot{H}_{x}^{-q}) }\leq \sum_{k=1}^{K}\alpha_k
\leq \frac{4\delta}{9} 4^K.
\end{aligned}
\end{equation*}

Furthermore, because $K$ is directly proportional to 
$l^{1/(1-2a)}\vee l^{1/(2-2(a+c))}$, then we have 
$$
\mathbb{E}\|\tilde{w}_\delta\|^2_{Z^{T_0\wedge \tilde{\iota}_{\delta,R}}_{b,0,q}\cap C([0,T_0\wedge \tilde{\iota}_{\delta,R}];L_x^2\cap \dot{H}_{x}^{-q}) }\leq \frac{4\delta}{9}C^{l^{1/(1-2a)}\vee l^{1/(2-2(a+c))}}.
$$

Therefore, for any $\varepsilon>0$ and $l>1$, there exists a
\begin{equation}\label{converging_speed}
\delta_{0}\lesssim \frac{\varepsilon}{l}\cdot C^{-\left(l^{1/(1-2a)}\vee l^{1/(2-2(a+c))}\right)},
\end{equation} 
such that
\begin{equation}\label{delta_convergence_result}
	\mathbb{E}\|\tilde{w}_\delta\|^2_{Z^{T_0\wedge \tilde{\iota}_{\delta,R}}_{b,0,q}\cap C([0,T_0\wedge \tilde{\iota}_{\delta,R}];L_x^2\cap \dot{H}_{x}^{-q}) }\leq\varepsilon/l,
\end{equation}
for any $\delta\in(0,\delta_0)$.

According to \eqref{uniform_moment_estimate_delta}, we know that $$\|w_\delta\|^2_{Z^{T_0}_{b,0,q}\cap C([0,T_0];L_x^2\cap \dot{H}_{x}^{-q})},\ \forall \delta\in(0,1)\ \text{and}\ \|w_0\|^2_{Z^{T_0}_{b,0,q}\cap C([0,T_0];L_x^2\cap \dot{H}_{x}^{-q})}$$
are uniformly integrable.

Precisely speaking, for any $A\in\mathcal{F}_{T_0}$, there exists a common constant $ C>0$ such that  
\begin{equation}\label{uniform_integrable}
\begin{aligned}
\mathbb{E}\|w_\delta\cdot\chi_{A}(\omega)\|^2_{Z^{T_0}_{b,0,q}\cap C([0,T_0];L_x^2\cap \dot{H}_{x}^{-q})} \leq& \left(\mathbb{E}\|w_\delta \|^4\right)^{1/2}\left(\mathbb{P}(A)\right)^{1/2}\\
\leq& C\cdot\left(\mathbb{P}(A)\right)^{1/2},\ \forall \delta\in[0,1).
\end{aligned}
\end{equation}

Note that 
\begin{equation}\label{partation_aux}
\begin{aligned}
\mathbb{E}\|\tilde{w}_\delta\|^2_{Z^{T_0}_{b,0,q}\cap C([0,T_0];L_x^2\cap \dot{H}_{x}^{-q})}
\lesssim&\mathbb{E}\|\tilde{w}_\delta\|^2_{Z^{T_0\wedge \tilde{\iota}_{\delta,\sqrt{l}R}}_{b,0,q}\cap C([0,T_0\wedge \tilde{\iota}_{\delta,\sqrt{l}R}];L_x^2\cap \dot{H}_{x}^{-q})}\\
&+\mathbb{E}\|w_\delta\|^2_{Z^{[T_0\wedge \tilde{\iota}_{\delta,\sqrt{l}R},T_0]}_{b,0,q}\cap C([T_0\wedge \tilde{\iota}_{\delta,\sqrt{l}R},T_0];L_x^2\cap \dot{H}_{x}^{-q})}\\
&+\mathbb{E}\|w\|^2_{Z^{[T_0\wedge \tilde{\iota}_{\delta,\sqrt{l}R},T_0]}_{b,0,q}\cap C([T_0\wedge \tilde{\iota}_{\delta,\sqrt{l}R},T_0];L_x^2\cap \dot{H}_{x}^{-q})}.
\end{aligned}
\end{equation}

Let us define 
$$
A_{\delta,R}=\left\{\omega:\|\tilde{w}_\delta\|_{Z^{T_0\wedge \tilde{\iota}_{\delta,R}}_{b,0,q}\cap C([0,T_0\wedge \tilde{\iota}_{\delta,R}];L_x^2\cap \dot{H}_{x}^{-q})}>R/3\right\},
$$
for any $R>0,\delta\in(0,1)$.

According to \eqref{delta_convergence_result} and \eqref{uniform_integrable}, it is obvious that there exists a constant $C>0$ such that 
\begin{equation}\label{auxi_A_R,delta}
\begin{gathered}
\mathbb{P}\left(A_{\delta,\sqrt{l}R}\right)\leq \frac{9\varepsilon}{l^2R^2},\\
\mathbb{E}\|w_\delta\cdot\chi_{A_{\delta,R}}\|^2_{Z^{T_0}_{b,0,q}\cap C([0,T_0];L_x^2\cap \dot{H}_{x}^{-q})}\leq C\frac{\varepsilon^{1/2}}{l},\\
\mathbb{E}\|w\cdot\chi_{A_{\delta,R}}\|^2_{Z^{T_0}_{b,0,q}\cap C([0,T_0];L_x^2\cap \dot{H}_{x}^{-q})}\leq C\frac{\varepsilon^{1/2}}{l},
\end{gathered}
\end{equation}
for any $\delta\in(0,\delta_0)$.

What's more, for any $\omega \in A^c_{\delta,\sqrt{l}R}$, we have 
\begin{equation}\label{distance}
	\tilde{\iota}_{\delta,R}\geq \iota_{0,2\sqrt{l}R/3} \vee \iota_{\delta,2\sqrt{l}R/3}.
\end{equation}

Therefore,  according to \eqref{partation_aux}, \eqref{distance}, \eqref{tail_moment_estimate}, \eqref{delta_convergence_result} and \eqref{auxi_A_R,delta}, we have 
\begin{equation}\label{strong_averaging_principle}
	\begin{aligned}
		&\mathbb{E}\|u-u_\delta\|^2_{Z^{T_0}_{b,0,q}\cap C([0,T_0];L_x^2\cap \dot{H}_{x}^{-q})}\\
		=&\mathbb{E}\|\tilde{w}_\delta\|^2_{Z^{T_0}_{b,0,q}\cap C([0,T_0];L_x^2\cap \dot{H}_{x}^{-q})}\\
		\lesssim&\mathbb{E}\|\tilde{w}_\delta\|^2_{Z^{T_0\wedge \tilde{\iota}_{\delta,\sqrt{l}R}}_{b,0,q}\cap C([0,T_0\wedge \tilde{\iota}_{\delta,\sqrt{l}R}];L_x^2\cap \dot{H}_{x}^{-q})}\\
		&+\mathbb{E}\|w_\delta\cdot\chi_{A_{\delta,\sqrt{l}R}}\|^2_{Z^{T_0}_{b,0,q}\cap C([0,T_0];L_x^2\cap \dot{H}_{x}^{-q})}\\
		&+\mathbb{E}\|w_\delta\cdot\chi_{A^c_{\delta,\sqrt{l}R}}\|^2_{Z^{[T_0\wedge \iota_{\delta,\sqrt{l}R/2},T_0]}_{b,0,q}\cap C([T_0\wedge \iota_{\delta,\sqrt{l}R/2},T_0];L_x^2\cap \dot{H}_{x}^{-q})}\\
		&+\mathbb{E}\|w\cdot\chi_{A_{\delta,\sqrt{l}R}}\|^2_{Z^{T_0}_{b,0,q}\cap C([0,T_0];L_x^2\cap \dot{H}_{x}^{-q})}\\
		&+\mathbb{E}\|w\cdot\chi_{A^c_{\delta,\sqrt{l}R}}\|^2_{Z^{[T_0\wedge \iota_{\sqrt{l}R/2},T_0]}_{b,0,q}\cap C([T_0\wedge \iota_{\sqrt{l}R/2},T_0];L_x^2\cap \dot{H}_{x}^{-q})}\\
		\leq& \frac{2}{l}(\varepsilon+2C\varepsilon^{1/2}),
	\end{aligned}
\end{equation}
for any $\delta\in(0,\delta_0)$.

Hence, we have finished the proof of Theorem \ref{thm:result_with_classical_workspace}.

\end{proof}

\section*{Acknowledgements}
Jie Chen was supported by the NSFC grants 12301116. Fan Gu  was supported by the NSFC grants 12601272.

\section*{Declaration of Generative AI Use}
The authors report generative Artificial Intelligence (AI) tools were not used in their research or preparation of this manuscript.

\phantomsection
\bibliographystyle{amsplain}
\addcontentsline{toc}{section}{References}
\bibliography{reference}

\scriptsize\textsc{Jie Chen: School of Sciences, Jimei University, Xiamen 361021, P.R. China}

\textit{E-mail address}: \textbf{jiechern@jmu.edu.cn}
\vspace{20pt}

\scriptsize\textsc{Fan Gu: School of Statistics and Mathematics, Central University of Finance and Economics,  Beijing 102206, P.R. China}

\textit{E-mail address}: \textbf{gufan@amss.ac.cn}
\vspace{20pt}

\scriptsize\textsc{Boling Guo: Institute of Applied Physics and Computational Mathematics, Beijing 100088, P.R. China}

\textit{E-mail address}: \textbf{gbl@mail.iapcm.ac.cn}
\vspace{20pt}

\scriptsize\textsc{Xianliang Zhao: Beijing International Center for Mathematical Research, Peking
University, Beijing 100871, P.R. China.}
%}

\textit{E-mail address}: \textbf{xzhaomath@gmail.com}
\vspace{20pt}

\end{document}